# Sun Bin's Legacy

Dana Mackenzie
scribe@danamackenzie.com

**Introduction**

Sun Bin was a legendary Chinese military strategist who lived more than 2000 years ago. Among other exploits, he is crediting with helping his patron, general Tian Ji, defeat the King of Qi in a match consisting of three horse races. If Tian Ji had simply raced his top horse against the King's top horse, his second against the King's second, and his third against the King's third, he would have lost all three races.

But Sun Bin had an idea. He told Tian Ji to race his *worst* horse against the King's best, his best against the King's second-best, and his second-best against the King's worst. In this way, Tian Ji won two out of three races.

Of course, Tian Ji was a little bit lucky. If we denote his horses by $A$, $B$, and $C$ and the king's horses by $a$, $b$, and $c$, and if we rank the horses by speed, they happened to fall in the order $a > A > b > B > c > C$. It's easy to see with 20/20 hindsight that Sun Bin's strategy works here, because $A > b$ and $B > c$. Had the horses been in a different order, say $a > b > c > A > B > C$, then Sun Bin's strategy would not have worked (but neither would any other strategy!).

A key point to realize, though, is that Sun Bin's strategy does not depend on knowing the relative speeds of *all six* horses in advance. We only need to know the rankings of each side's horses: that is, we only need to know that $A > B > C$ and $a > b > c$. Given only this information, Sun Bin's strategy of racing $C$ against $a$, $A$ against $b$, and $B$ against $c$ is the optimal strategy in two different ways: It gives him the best odds of winning the match, and gives him the largest number of expected races won. I will justify these claims in this paper.

Now, let's bring Sun Bin's problem into the twenty-first century. What happens if we have a match of $N$ horses against $N$ horses? Can we find an optimal strategy for winning the match? How about for winning the largest expected number of match races? The answer to both of these questions is yes, and the goal of this paper is to answer both questions. The second one is much harder, and (in my opinion) much more interesting mathematically.

To define the problem a little better, I prefer to phrase it in terms of a card game. The deck has $2N$ cards, with face values from 1 to $2N$ ($2N$ being high). Player 1 receives $N$ cards, face down, and player 2 likewise. Neither is allowed to see the cards in his hand or the other player's hand. However, the cards are placed in rank order in front of each player, so that each player knows the relative ranking of his own cards and the opponent's cards. On the first trick, player 1 plays one card (still face down) and player 2 plays one against it. Play continues in this fashion (player 1 always going first!) until all the cards have been played, and then both players reveal their cards. The winner is the one who takes the most tricks. What is player 2's optimal strategy?

It is easy to see that this is equivalent to an "idealized" version of the horse problem, in which the faster horse always wins the race. I prefer the card version because (as is well known) horses do not always race according to form, and faster horses sometimes lose to slower horses. In the card version, there are no such ambiguities: card 7 always beats card 6, and that is that.

To me, the *N*-card (or *N*-horse) version of Sun Bin's problem is extremely natural, and it is a bit of a mystery why it seems to be nearly absent from the mathematical literature. The only reference I have been able to find is [AGY], written in 1979, and even that reference is cursory. The authors simply noted that the problem reduces to a linear assignment problem, and therefore there exist fast computer algorithms (the Hungarian Algorithm) to solve it.

A college friend of mine, Howard Stern, independently posed this problem in his first year of graduate school, in 1980. He made, in my opinion, some extremely impressive progress towards a solution, and arrived at a correct conjecture for the general strategy, but he was unable to prove it. In the three decades since then, he has showed the problem to a number of mathematicians and computer scientists, always thinking that somebody would know a solution or a general theorem that would solve the problem. However, none of them did. Finally, he asked me in 2012, and very soon I was hooked on Stern's problem! The proof of his conjecture turned out to be a fascinating mix of group theory, probability, and combinatorics.

What's more, I believe that the solution has some relevance to the general theory of linear assignment problems. Stern's original, unpublished work from 1980 appears to be a novel observation about what I call mixed-Monge optimization problems. And the first step of my proof of Howard's conjecture also involves a general result about symmetric mixed-Monge optimization problems.

Gary Antonick wrote a post about Stern's problem in his "Numberplay" blog for the *New York Times* on January 13, 2014. It became his most-commented-on blog post in more than a year. I am indebted to one of his readers, a reader known to me only as Lee from London, who pointed out the ancient Chinese legend of Sun Bin. This is surely the first appearance of the problem in recorded history, so I think it is only appropriate to call it "Sun Bin's Legacy."

At this point I would strongly encourage readers to play the card game and see if they can figure out the strategy themselves, before reading on. Antonick's post [A] includes a wonderful applet by Gary Hewitt that will enable you to play against the computer online with 3 to 7 cards.

The main theorem of this paper is as follows:

**Main Theorem.** For sufficiently large $N$, the optimal strategy for player 2 is to play his cards in the order $(N, N-1, \ldots, N-k+1, 1, 2, \ldots, k)$ for some $k$. In other words, he plays his $N$-th card against player 1's first card, his $(N-1)$-th card against player 1's second, etc. Note that $k$ represents the number of tricks that player 2 should (try to) "throw" or lose on purpose. This strategy is optimal in the sense of maximizing the *expected number of tricks won*. The optimal number $k = k^*(N)$ is given by the following formula:

$$k^*(N) = \sup\left\{k : \binom{N}{k-1}^2 + \sum_{j=0}^{N-k}\binom{2N}{j} \geq \binom{2N}{N}\right\}.$$

**Comments on the Main Theorem:**

(1) There is also a somewhat simpler "approximate" formula $k \approx k(N)$:
$$k(N) = \sup\left\{k : \sum_{j=0}^{N-k} \binom{2N}{j} \geq \binom{2N}{N}\right\}.$$
It is approximate in the sense that $k^*(N) - k(N) = 0$ or $1$. In fact, I do not know a single value of $N$ for which $k(N) \neq k^*(N)$.

(2) At present, my best rigorous estimate for "sufficiently large $N$" is for $N > 10{,}000{,}000$. I believe this estimate can be improved to about $N > 400$, using more careful versions of the estimates in this paper. The theorem is almost certainly true for all $N \leq 400$ as well. Using the computer, Stern has verified it for all $N \leq 60$. If somebody with a supercomputer can verify it for $60 < N \leq 400$, then I would cheerfully undertake the work of formalizing the proof for $N > 400$.

(3) Stern conjectured the general form of the optimal strategy in 1980, but did not make a conjecture for the optimal number $k^*(N)$ of tricks to throw. At that point there was not enough data to make a conjecture, and it is highly unlikely that anyone would have come up with the formula above anyway. It was a complete shock to me that I was able to derive an exact formula.

Here is a table of the optimal number of tricks to throw for small values of $N$:

| N | k | N | k |
|---|---|---|---|
| 2 | n/a | 9 | 2 |
| 3 | 1 | 10 | 2 |
| 4 | 1 | 11 | 2 |
| 5 | 1 | 12 | 2 |
| 6 | 1 | 13 | 3 |
| 7 | 2 | 14 | 3 |
| 8 | 2 | 15 | 3 |

And here are two "worked examples," showing the first two jumps in $k^*(N)$.

**Example 1:** $N = 7$. Here the approximate formula tells us to look up the 14-th row of Pascal's triangle and add the terms until we get a sum that is greater than the central element. We find that
$$1 + 14 + 91 + 364 + 1001 + 2002 = 3473 > 3432.$$
Then the approximate number of tricks to throw is $(N+1)$ minus the number of terms added: in this case $(7 + 1) - 6 = 2$.

The exact computation simply involves adding one term to the above inequality. We note that
$$1 + 14 + 91 + 364 + 1001 + 2002 + (7)^2 = 3522 > 3432.$$
Thus the exact number of tricks to throw is at least 2. On the other hand, if we try throwing one more, we get
$$1 + 14 + 91 + 364 + 1001 + (21)^2 = 1912 < 3432.$$
Thus the exact number of tricks to throw is at most 2, and hence the exact formula agrees with the approximate formula.

**Example 2:** $N = 13$. Now we look up the 26[th] row of Pascal's triangle and add up the terms until we get a sum that is greater than the center element. We find that

$$1 + 26 + 325 + 2600 + 14950 + 65780 + 230230 + 657800 + 1562275 + 3124550 + 5311735 = 10970722 > 10400600.$$
Therefore the approximate number of tricks to throw is ($N$ + 1) minus the number of terms added, i.e. (13 + 1) – 11 = 3.

The exact computation involves adding one term from the 13$^{th}$ row, squared. Because
$$1 + 26 + 325 + 2600 + 14950 + 65780 + 230230 + 657800 + 1562275 + 3124550 + 5311735 + (78)^2 > 10400600$$
we can be certain that the exact number of tricks to throw is at least 3. And because
$$1 + 26 + 325 + 2600 + 14950 + 65780 + 230230 + 657800 + 1562275 + 3124550 + (286)^2 < 10400600$$
the exact number of tricks to throw is less than 4. Hence the exact number of tricks to throw is 3, which agrees with the approximate computation.

These two examples are completely typical. The additional "nuisance term" from the $N$-th row of Pascal's triangle, even though it is squared, is dwarfed by the largest terms from the ($2N$)-th row. This is why the "approximate" formula agrees with the exact formula in every case I know of.

(4) It is also of interest to derive upper and lower bounds for the exact number of tricks to throw. After all, if you are playing the game with $N$ = 200 cards, it may not be so easy to look up the 400$^{th}$ row of Pascal's triangle! I prove the following estimate in this paper:

**Theorem:** If $N$ > 400, and if the optimal permutation is of the type in the Main Theorem, then the optimal number of tricks to throw satisfies the inequalities
$$\sqrt{N \ln N / 4} < k^*(N) < \sqrt{N \ln N / 2}.$$
Computer calculations by Stern show that these inequalities hold for 400 ≥ $N$ ≥ 91 as well. As $N \to \infty$, $k^*(N) \sim \sqrt{N \ln N / 2}$. It is interesting that this asymptotic limit is approached extremely slowly. Stern's computer calculations show that all the way up to $N$ = 500, the ratio $k^*(N)/\sqrt{N \ln N}$ is closer to 0.5 than it is to 0.7071..., its eventual limit. (For $N$ = 500, the ratio is 0.559...)

Note that the hypothesis "if the optimal permutation is of the type in the Main Theorem" is added here only because I do not have a proof of the Main Theorem yet for $N$ less than 10 million. The asymptotic result as $N \to \infty$ is true unconditionally.

Finally, as mentioned briefly above, there is a second version of Sun Bin's Legacy, which is to find the strategy that guarantees the highest probability of winning a majority of tricks, regardless of the number of tricks won. Curiously, neither Stern nor I worked seriously on this question. In my case, this was because I expected the majority-of-tricks problem to be harder, because the objective function is nonlinear.

Imagine my astonishment when, within one day of Gary Antonick's post going up on the "Numberplay" blog, one of his readers found the optimal strategy for the majority-of-tricks problem! Here I assume $N$ = 2$n$+1 is odd. Then reader Bill Courtney showed that the optimal strategy is to throw $n$ tricks. Thus player 2 pairs his top ($n$+1) cards against player 1's bottom ($n$+1) cards, in order. It is easy to see that if there is *any way at all* to win ($n$+1) tricks, then this strategy will do so. The proof is left to the reader (or see Courtney's comment to [A]).

While Courtney's strategy maximizes the probability of winning a simple majority, it is extravagantly wasteful on the level of tricks. It will on average lose nearly half the tricks. By contrast, the "Pascal's triangle" strategy described above will on average lose only about

$\sqrt{N \ln N / 2}$ tricks. By playing with a large enough deck, you can win as close as you want to 100 percent of the tricks!

The outline of the rest of the paper is as follows:

**I. Basic Results, Mixed Monge Matrices and the Shape Theorem.**

This section sets up the problem as a linear assignment problem, shows that the objective function is given by a "mixed Monge matrix," and derives a weak form of the optimal strategy. In particular, I show that the optimal strategy always involves throwing some tricks in reverse order, and playing the rest of the tricks in normal order. However, there may be "gaps" in the thrown tricks. Most of the work in this section is due to Stern (unpublished).

**II. The Symmetry Lemma.**

The objective function in section I leads to a mixed Monge matrix that is symmetric about the "anti-main diagonal," and skew-symmetric (after subtracting a constant from each entry) about the main diagonal. I exploit this symmetry to prove that *if* you have decided which tricks to throw (say tricks 1, 3, and 7) then the optimal strategy for these tricks is to play your *i*-th worst card against your opponent's *i*-th best card. Still, there may be gaps in the thrown tricks.

**III. The No-Gaps Theorem.**

In this section, which is the most technical one, I show that if *N* is large enough (at least 10 million) then the optimal strategy has no gaps. That is, you should throw tricks 1, 2, …, *k* for some *k*. Although I do not derive the best strategy in this section, the proof depends on knowing that a *very good* strategy is to sacrifice the first $\sqrt{N \ln N / 2}$ tricks. Roughly speaking, this strategy beats any strategy with gaps in it.

**IV. The Number of Tricks to Throw.**

In the last section, I explain the wonderful and totally unexpected connection between the expected number of tricks won and the 2*N*-th row of Pascal's triangle. I derive the exact number of tricks to throw, $k^*(N)$, the approximate number $k(N)$, and the asymptotic limit for both of them.

**Bibliography.**

[A] G. Antonick, "Stern-Mackenzie One-Round War," *New York Times* ("Numberplay" blog), Jan. 13, 2014, at http://wordplay.blogs.nytimes.com/2014/01/13/war/.
[AGY] A. Assad, Golden, B., and Yee, J. Scheduling players in team competitions. *Proceedings of the Second International Conference on Mathematical Modelling (St. Louis, MO, 1979),* Vol I, pp. 369-379. Rolla, MO: Univ. of Missouri-Rolla, 1980.

# I. Basic Results, Mixed Monge Matrices and the Shape Theorem

I will begin with a review the rules of the card game, which Gary Antonick named Stern-Mackenzie One Round War.

The game requires three people, one of whom serves only as the dealer. The dealer has a deck of $2N$ cards, numbered one through $2N$. She shuffles the cards and then deals $N$ cards to each player, making sure, as she deals the cards, to place them in order from highest to lowest (she can look at the cards as she's dealing them without the players seeing their values).

Cards are then put face-down into the middle of the table. Player 1 plays a card, without turning it over; then player 2 plays a card, without turning it over. Those two cards constitute the first trick. Play continues in the same way until all the cards have been paired up. Only then do the players turn over all the cards. The higher card takes each trick. The player who wins the most tricks wins the game.

This paper is about answering the question: What is the optimal strategy for player 2, in the sense of maximizing the expected number of tricks won? (As discussed in the introduction, this is not the same as the optimal strategy for maximizing the probability of winning the game.)

A "strategy" is simply a permutation $\pi$ of the numbers 1 through $N$. If we rank player 2's cards from bottom to top, with 1 being his worst card and $N$ being his best, then strategy $\pi$ for player 2 is to play his card $i$ against player 1's card $\pi(i)$. Thus, for example, if $N = 3$, then Sun Bin's strategy corresponds to the permutation $\pi$ such that $\pi(1) = 3$, $\pi(2) = 1$, and $\pi(3) = 2$. In group theory, this permutation is written (1 3 2).

[*Note:* The rank $i$ of an individual card is called the "order statistic" of that card, and should not be confused with the actual face value of the card, which is a number between 1 and $2N$.]

How do we measure the expected number of tricks won by a given strategy? Following Stern's approach, let $p_{ij}$ denote the probability that the $i$-th worst card in player 2's hand beats the $j$-th worst card in player 1's hand. Then the expected number of tricks won by the $i$-th card, using strategy $\pi$, is simply $p_{i\pi(i)}$, and the total expected number of tricks won is

$$\sum_{i=1}^{N} p_{i\pi(i)}.$$

A permutation $\pi$ can be associated with a matrix of 0's and 1's, which I will denote by $M_\pi$, by putting a "1" in the $ij$-th position if $\pi(i) = j$ and a 0 otherwise. If we let **P** denote the matrix whose $ij$-th entry is $p_{ij}$, then the total expected number of tricks won by strategy $\pi$ is $M_\pi \cdot \mathbf{P}$.

The problem of finding the permutation $\pi$ that maximizes this inner product is called a "linear assignment problem." (See [B].) It is a very famous type of problem in operations research, which from an algorithmic point of view is well understood. The Hungarian Algorithm solves the problem in $O(N^4)$ steps, and it has been modified to work in $O(N^3)$ steps.

I would like to comment on two related papers to highlight ways in which they differ from this one. In the Introduction I mentioned [AGY], the only paper I have found that explicitly considers something more or less like Sun Bin's problem. Assad et. al. derive the above formula for the expected number of tricks won, point out that it is a linear assignment problem, which has efficient solution algorithms, and say little more about it. One reason is that they are operating with different assumptions. Their context is a sporting match, in which the coach of team A (the equivalent of our player 2) has a perfect scouting report on the other team. Thus he has empirical knowledge of the numbers $p_{ij}$, the probability that his $i$-th player will defeat the other team's $j$-th player, for each $i$ and $j$. In contrast, we do not assume any prior knowledge of the opponent's cards. Our matrix **P** is, in some sense, an uninformative prior for those probabilities. These are the probabilities that the scouts would report back to the coach if they didn't actually watch the other team's players but just spent the afternoon at a bar.

Because Assad is looking at a very general probability matrix **P**, it is natural that he would not be able to offer any better solution than the Hungarian Algorithm or any other general-purpose algorithm for solving linear assignment problems. However, because we will look at a specific **P** (or more precisely, a specific $\mathbf{P}^N$ for any given number of cards $N$), it is reasonable to ask for more than a solution algorithm—we ask for an actual solution.

Another contrasting paper is [DL]. Again this paper considers a match between two teams, or more generally two populations from which randomly chosen groups of $N$ objects are compared in pairs. David and Liu are looking for *fair matchings*. These are strategies $\pi$ that ensure that the expected number of paired comparisons won by each population is $N/2$, assuming that the distribution of strengths in each population is the same. One example of a fair matching is the identity permutation, which matches the strongest player against the strongest, the second against the second, etc. By contrast, in this paper I am looking for the *most unfair* matching.

What do the matrices $\mathbf{P}^N$ look like? Here are the first four, starting with $N = 2$:

$$P^2 = \begin{bmatrix} 3 & 1 \\ 5 & 3 \end{bmatrix}$$

$$P^3 = \begin{bmatrix} 10 & 4 & 1 \\ 16 & 10 & 4 \\ 19 & 16 & 10 \end{bmatrix}$$

$$P^4 = \begin{bmatrix} 35 & 15 & 5 & 1 \\ 55 & 35 & 17 & 5 \\ 65 & 53 & 35 & 15 \\ 69 & 65 & 55 & 35 \end{bmatrix}$$

$$P^5 = \begin{bmatrix} 126 & 56 & 21 & 6 & 1 \\ 196 & 126 & 66 & 26 & 6 \\ 231 & 186 & 126 & 66 & 21 \\ 246 & 226 & 186 & 126 & 56 \\ 251 & 246 & 231 & 196 & 126 \end{bmatrix}, \text{etc.}$$

For any readers who might be interested, matrices $\mathbf{P}^6$ through $\mathbf{P}^8$ are given in Table 1.1.

For reasons of convenience I have scaled up each of the matrices by a factor of $\binom{2N}{N}$, the total number of possible hands in Stern-Mackenzie One-Round War. Thus, for instance, $p_{24}^5 = 26$ means that there are 26 possible hands in five-card One-Round War in which player 2's second-worst card beats player 1's fourth-worst card. The entry in the northeast corner is always 1, because there is only one hand in which player 2's worst card beats player 1's best card. That will occur if and only if player 2 receives the cards with face values ($N$+1) through $2N$.

Clearly the re-scaling by $\binom{2N}{N}$ does not affect Sun Bin's problem in any material way. The goal of the problem is to circle one entry in each row and column of **P** in such a way that the sum of the circled elements is maximized.

**Lemma 1.1.** (Stern) $p_{ij}^N = \sum_{k=j}^{N} \binom{k+i-1}{i-1}\binom{2N-k-i}{N-i}$.

As much as possible we will try to avoid using this formula, but it is sometimes inevitable. Even more confusingly, there is a *second* formula whose proof we will postpone until Section 4:

**Lemma 1.2 = Lemma 4.3.** (Mackenzie) $p_{ij}^N = \sum_{k=1}^{j} \binom{i+j-1}{k-1}\binom{2N-i-j+1}{N-k+1}$. A significant special case occurs when $i+j = N+1$. In this case I will frequently say that the entry $p_{ij}$ lies on the "anti-main diagonal" (AMD). If $p_{ij}$ lies on the AMD with $i<j$, then $p_{ij} = 1 + \binom{N}{1}^2 + ... + \binom{N}{i-1}^2$.

While they look quite similar, Lemmas 1.1 and 1.2 are actually different. Stern's formula "telescopes" along the rows, which means that $p_{ij} - p_{i,j+1}$ has a simple expression as a product of two binomial coefficients. My formula telescopes along the diagonals $i+j$ = constant, which means that $p_{ij} - p_{i-1,j+1}$ can also be expressed as a product of two binomial coefficients.

**Lemma 1.3.** (Stern) $p_{ij} + p_{ji} = \binom{2N}{N}$.

Thus, if we subtract $\frac{1}{2}\binom{2N}{N}$ from every entry of **P**, we get a skew-symmetric matrix. It is really a matter of taste whether to do this or not, but I choose not to do it in order to keep the formulas for $p_{ij}$ simpler.

**Lemma 1.4.** (Stern) $p_{ij} = p_{N-j+1,N-i+1}$.

**Proof.** The left-hand side is the probability that player 2's $i$-th worst card beats player 1's $j$-th worst. Because the probabilities don't depend on the labeling of players, this is the same as the probability that player 1's $i$-th worst card beats player 2's $j$-th worst card. Because the probabilities are invariant under simultaneous interchange of the English words "best" ⇔ "worst" and "beats" ⇔ "loses to", this is the same as the probability that player 1's $i$-th best card loses to player 2's $j$-th best card. Finally, the case just described is identical to player 2's ($N$+1-$j$)-th worst card beating player 1's ($N$+1-$i$)-th worst card. □

Lemma 1.4 implies that the matrix **P** is symmetric about the anti-main diagonal, which leads from northeast to southwest. This symmetry property will be of the utmost importance in section 2 of this paper, when it will be used to show that the solution $M_\pi$ of Sun Bin's problem is likewise symmetric about the anti-main diagonal.

**Lemma 1.5.** (Stern) The entries in $\mathbf{P}^N$ increase from top to bottom along the columns, and from right to left along the rows.

This is pretty obvious; it means that higher ranked cards in player 2's (resp. player 1's) hand are more likely to beat any given card in the opponent's hand.

Although the formulas for individual elements of $\mathbf{P}$ are relatively messy, certain linear combinations of those elements are very nice. For example, the sums of the columns increase linearly from right to left:

**Lemma 1.6.** (Mackenzie) The sum of the elements in the ($N$-$j$+1)-th column or the $j$-th row of $\mathbf{P}^N$ is $j\binom{2N}{N-1}$.

**Proof.** For $j = 1$, we note that the elements in the $N$-th column are $\binom{N}{0}, \binom{N+1}{1}, \ldots, \binom{2N-1}{N-1}$, and it is a well-known identity in Pascal's triangle that the sum of these is $\binom{2N}{N-1}$. To complete the proof, we need to show that the difference between any two adjacent column sums is also $\binom{2N}{N-1}$. Using Stern's Lemma 1.1, this difference is equal to
$$\sum_{i=1}^{N} \binom{j+i-1}{i-1}\binom{2N-i-j}{N-i}.$$
To simplify this, use the general identity for binomial coefficients:
$$\sum_{i=0}^{m} \binom{i+j}{i}\binom{m-i+k}{m-i} = \binom{m+j+k+1}{m}.$$
This identity probably cannot be called "well-known," but it is easy to prove it. The right-hand side counts the number of ways to choose $m$ numbers from the set $\{1, 2, \ldots, m+j+k+1\}$. We ask, "Where is the ($j$+1)-th missing number?" If the ($j$+1)-th missing number is at position ($i$+$j$+1), that means we have chosen $i$ numbers from the set $\{1, 2, \ldots, i+j\}$. The first term in parentheses counts the number of ways to do this. There are then $k$ missing numbers in the remaining positions, $\{i+j+2, \ldots, m+j+k+1\}$, and the number of ways to choose these is given by the second term in parentheses. Adding up over all of the possible locations for the ($j$+1)-th missing number gives the summation. Substituting $m = N$-1 and $k = N$-$j$ produces the desired formula for the difference between column sums. Finally, the statement for rows follows from the statement for columns by symmetry (Lemma 1.4). $\square$

A very important quantity for Sun Bin's problem is the sum of the entries along any diagonal parallel to the main diagonal. Let $s_{kN}$ denote the sum of all the elements of $\mathbf{P}^N$ that lie $k$ steps above the main diagonal, and let $\bar{s}_{kN}$ denote the sum of all the elements that lie $k$ steps below the main diagonal.

**Lemma 1.7.** (Mackenzie) $s_{kN} = \dfrac{N}{2}\binom{2N}{N} - k\binom{2N}{0} - k\binom{2N}{1} - \ldots - k\binom{2N}{N-k} - (k-1)\binom{2N}{N-k+1} - \ldots - \binom{2N}{N-1}$

and $\bar{s}_{kN} = \left(\dfrac{N}{2} - k\right)\binom{2N}{N} + k\binom{2N}{0} + k\binom{2N}{1} + \ldots + k\binom{2N}{N-k} + (k-1)\binom{2N}{N-k+1} + \ldots + \binom{2N}{N-1}$.

The proof of this lemma is much more difficult and will be given in Section 4. However, it is very easy to check for the first few **P** matrices. For example, for $N = 4$ we have:

$s_{04} = 2 \cdot 70 = 140$, $s_{14} = 2 \cdot 70 - 1 - 8 - 28 - 56 = 47$,
$s_{24} = 2 \cdot 70 - 2 \cdot 1 - 2 \cdot 8 - 2 \cdot 28 - 56 = 10$, and $s_{34} = 2 \cdot 70 - 3 \cdot 1 - 3 \cdot 8 - 2 \cdot 28 - 56 = 1$.
All of these diagonal sums can easily be confirmed by referring to $\mathbf{P}^4$, given above.

The next lemma is Stern's most remarkable discovery about the **P** matrices.

**Lemma 1.8. (The Mixed Monge Property.)** (Stern) Consider any four entries of **P** that form the vertices of a rectangle. Call them "northeast" (NE), "northwest" (NW), "southeast" (SE) or "southwest" (SW) according to their locations. If three or more of the four entries lie on or above the main diagonal, then NE + SW ≥ NW + SE. If three or more of the four entries lie on or below the main diagonal, then NW + SE ≥ NE + SW. If two entries lie on the main diagonal, then NW + SE = NE + SW. The test is indeterminate if two entries lie above the main diagonal and two lie below.

**Proof.** First we prove Lemma 1.6 for four adjacent entries that form a unit square. Thus we have to compute $p_{ij} - p_{i+1,j} - p_{i,j+1} + p_{i+1,j+1}$. Lemma 1.1 is ideal for this purpose, because it telescopes along the rows. The details are left to the reader, but the result is as follows:

$$p_{ij} - p_{i+1,j} - p_{i,j+1} + p_{i+1,j+1} = \binom{2N-i-j}{N-j}\binom{i+j}{j}\frac{N(i-j)}{(2N-i-j)(j+1)}.$$

As messy as it appears, only one thing matters about this formula: everything on the right-hand side is positive except for $i - j$. If $i > j$, which is the case where three or more of the vertices lie below the main diagonal, then we conclude that NW + SE − NE − SW > 0. If $i < j$, which is the case when three or more of the vertices lie above the main diagonal, we conclude that NW + SE − NE − SW < 0. If $i = j$, which is the case when two of the vertices lie on the main diagonal, we have NW + SE = NE + SW.

Now, given four arbitrary entries in **P** forming a rectangle, we partition the rectangle into unit squares. In all of the determinate cases, we can use the skew symmetry of **P** to cancel out submatrices that lie above the main diagonal with symmetrically placed submatrices below the main diagonal. The remaining unit squares after the cancellation will either lie all above the main diagonal or all below the main diagonal, and the sums NW + SE − NE − SW will all be negative or all positive correspondingly. □

**Comments on Lemma 1.8.**

**(1)** The reason that the proof does not work for the indeterminate case is that after cancellation of symmetrically placed unit squares, there will be some uncancelled squares both above and below the main diagonal. The solution of Sun Bin's problem would be vastly simpler if there were some universal way to deal with elongated rectangles that cross the main diagonal (see Figure 1.1). Stern tried and failed to find a rule, and so did I.

**(2)** The "Monge property" is a known concept in linear assignment problems (see [B]). A Monge matrix is one in which the rule NW + SE − NE − SW > 0 holds universally. In this case, the optimal permutation follows the main diagonal, i.e. it is the identity permutation. The proof is simple. Given a permutation matrix that isn't the diagonal matrix, it must have two nonzero elements that are the SW and NE corners of a rectangle. That is, $i < j$ but $\pi(i) > \pi(j)$. Now just switch them, i.e., define a new permutation $\pi'$ such that $\pi'(i) = \pi(j)$ and $\pi'(j) = \pi(i)$. Now the matrix $M_{\pi'}$ is identical to $M_\pi$, except that the old SW and NE corners have been replaced by NW and SE corners. According to the Monge property, this substitution increases the value of the objective function. Eventually we must run out of pairs to switch because the objective function is bounded. When the algorithm terminates, it results in a permutation matrix with no SW-NE corners, or all NW-SE corners. The only such permutation matrix is the identity.

I call Stern's matrix **P** a "mixed Monge matrix" because it can't make up its mind. Below the diagonal, it has the Monge property. Above the main diagonal, it has the anti-Monge property. Stern's most important result for Sun Bin's problem was to note that the assignment problem for such matrices likewise has a solution with a canonical form. However, the form is not unique, as it would be for a pure Monge matrix.

**Theorem 1.9. (The Shape Theorem.)** (Stern) If **P** is an $N \times N$ mixed Monge matrix and $N > 2$, then any permutation matrix $M_\pi$ that maximizes the objective function $M_\pi \cdot \mathbf{P}$ must satisfy the following conditions:
(1) No nonzero elements on the main diagonal;
(2) All the nonzero elements below the main diagonal form a decreasing sequence (i.e., they run from NW to SE, but they do not have to lie on a straight line).
(3) All the nonzero elements above the main diagonal form an increasing sequence (i.e., they run from SW to NE, but do not have to lie on a straight line).

**Proof.** Starting with any candidate permutation matrix, first look for the left-most nonzero entry on the main diagonal. If it is not in the first column, there is a nonzero entry in the column immediately to its left. If that entry is above the main diagonal, switching NW + SE $\Rightarrow$ NE + SW will increase the objective function by Lemma 1.8. If the entry in the adjacent column is below the main diagonal, switching NE + SW $\Rightarrow$ NW + SE will increase the objective function.

If the left-most nonzero entry on the main diagonal is in the (1, 1) position, look for the first row that has an entry not on the main diagonal. It must in fact lie above the main diagonal. In that case, it lies SE of the (1, 1) position and switching the two entries will again increase the objective function by Lemma 1.8.

If the left-most nonzero entry on the main diagonal is in the (1, 1) position and every entry below it is on the main diagonal, then the permutation is just the identity permutation. In that case, switch any pair from NW + SE to NE + SW. This will not decrease the objective function, and it will lead to the situation of step one or two, in which the objective function can be increased.

Thus any matrix that has a left-most entry on the main diagonal is not optimal. The contrapositive statement is: Any optimal matrix cannot have a left-most entry on the main diagonal. In other words, it cannot have *any* entry on the main diagonal, which proves (1).

Now consider any candidate matrix with no entries on the main diagonal. Look either for NW + SE pairs above the main diagonal, or NE + SW pairs below the main diagonal. If you can find such a pair, make a substitution NW + SE $\Rightarrow$ NE + SW above the main diagonal or make a substitution NE + SW ⇒ NW + SE below the main diagonal. In either case, you will strictly increase the objective function, so the candidate matrix could not have been optimal. If you cannot find such a pair, the permutation matrix has the form described in the theorem. □

**Comments on Theorem 1.9**

**(1)** Notice that the proof of Theorem 1.9 is completely independent of the specific values of the **P** matrix. As far as I know, no one has ever observed this property of mixed Monge matrices before. It could be called "Stern's Theorem."

**(2)** If we interpret Theorem 1.9 in the context of Sun Bin's problem, property (1) means that you should never play a card against an equally-ranked card of your opponent. Thus your strategy consists of dividing your cards into two sets: a set $T$ of cards that will be "thrown" by pairing them against higher-ranked cards in your opponent's hand, and a set $S$ of cards that will be paired against lower-ranked cards. (Perhaps we could call these the "straight" cards.) Property (2) says that set $T$ forms a SW-to-NE sequence in the permutation matrix, which means they are paired against your opponent's cards in *reverse* order. Property (3) says that set $S$ forms a NW-to-SE sequence in the permutation matrix, which means that these should be paired against your opponent's cards in *proper* order.

To give an example, if $N = 5$ there are four permutations that satisfy the conclusions of Theorem 1.7. The optimal one is the permutation (1 5 4 3 2), in which you sacrifice only one trick. You play your card 1 (your weakest card) against the opponent's card 5 (his strongest), your card 5 against his card 4, your card 4 against his card 3, and so on. The set $T$ of thrown cards is {1} and the set $S$ of straight cards is {2, 3, 4, 5}.

The three other permutations that satisfy the conclusions of Theorem 1.9 but are not optimal are (1 5 3)(2 4), (1 5 4 2 3) and (1 5 3 4 2). The last one is particularly interesting because the set $T$ of thrown cards has a gap: $T = \{1, 3\}$. (Your card 1 is sacrificed against the opponent's card 5, and your card 3 is sacrificed against the opponent's card 4.) In Section 3 we will show that an optimal permutation cannot have gaps. Actually, the last two permutations will also be ruled out by Section 2, because they are not symmetric about the anti-main diagonal.

**(3)** A corollary of Theorem 1.9 is that you should always play your worst card against the opponent's best; i.e., in an optimal permutation, $\pi(1) = N$.

**(4)** Although I will not need this result anywhere, it is an interesting problem to show that the number of permutation matrices that satisfies the conclusion of Theorem 1.9 is $2^{N-3}$. (For example, when $N = 5$ there are four possibilities.) Stern made this conjecture around 1980 and it was proven by Evan O'Dorney in 2012. Thus we have reduced the number of candidate strategies from $N!$ to $2^{N-3}$, a quite dramatic reduction! Still, there is a long way to go to reduce the number to 1.

$N = 5$:

| 126 | 56  | 21  | 6   | 1   |
|-----|-----|-----|-----|-----|
| 196 | 126 | 66  | 26  | 6   |
| 231 | 186 | 126 | 66  | 21  |
| 246 | 226 | 186 | 126 | 56  |
| 251 | 246 | 231 | 196 | 126 |

$(21 + 66 + 126 + 66 + 21) - (6 + 26 + 6) = 262$

$N = 6$:

| 462 | 210 | 84  | 28  | 7   | 1   |
|-----|-----|-----|-----|-----|-----|
| 714 | 462 | 252 | 112 | 37  | 7   |
| 840 | 672 | 462 | 262 | 112 | 28  |
| 896 | 812 | 662 | 462 | 252 | 84  |
| 917 | 887 | 812 | 672 | 462 | 210 |
| 923 | 917 | 896 | 840 | 714 | 462 |

$N = 7$:

| 1716 | 792  | 330  | 120  | 36   | 8    | 1    |
|------|------|------|------|------|------|------|
| 2640 | 1716 | 960  | 456  | 176  | 50   | 8    |
| 3102 | 2472 | 1716 | 1016 | 491  | 176  | 36   |
| 3312 | 2976 | 2416 | 1716 | 1016 | 456  | 120  |
| 3396 | 3256 | 2941 | 2416 | 1716 | 960  | 330  |
| 3424 | 3382 | 3256 | 2976 | 2472 | 1716 | 792  |
| 3431 | 3424 | 3396 | 3312 | 3102 | 2640 | 1716 |

$(792 + 1716 + 2472 + 1716 + 1016 + 491 + 176 + 36) - (330 + 960 + 456 + 176 + 50 + 8) = 6435$

$N = 8$:

| 6435  | 3003  | 1287  | 495   | 165   | 45    | 9    | 1    |
|-------|-------|-------|-------|-------|-------|------|------|
| 9867  | 6435  | 3663  | 1815  | 765   | 261   | 65   | 9    |
| 11583 | 9207  | 6435  | 3915  | 2025  | 849   | 261  | 45   |
| 12375 | 11055 | 8955  | 6435  | 3985  | 2025  | 765  | 165  |
| 12705 | 12105 | 10845 | 8885  | 6435  | 3915  | 1815 | 495  |
| 12825 | 12609 | 12021 | 10845 | 8955  | 6435  | 3663 | 1287 |
| 12861 | 12805 | 12609 | 12105 | 11055 | 9207  | 6435 | 3003 |
| 12869 | 12861 | 12825 | 12705 | 12375 | 11583 | 9867 | 6425 |

**Table 1.1.** The matrices $\mathbf{P}^5$ through $\mathbf{P}^8$. Also illustrated is the "hook-sum" rule for recursively generating $\mathbf{P}^{N+1}$ from $\mathbf{P}^N$ (see Theorem 4.2).

## II. The Symmetry Lemma.

The Shape Theorem reduces the number of permutations that could be maximal, from $N!$ to $2^{N-3}$. The next step will be the Symmetry Lemma, which states that the matrix representation of an optimal permutation must be symmetric about the Anti-Main Diagonal (AMD). This will reduce the number of possible maximizers to about $2^{N/2}$, and it will give them a very specific structure.

Before proceeding with the proof, let us make some philosophical observations. First, because the matrix **P** is symmetric about the AMD, the objective function $F(\pi) = M_\pi \cdot \mathbf{P}$ is also symmetric. That is, if $\pi'$ denotes the permutation whose matrix is the same as $M_\pi$ reflected through the AMD, then $F(\pi) = F(\pi')$. Thus it is natural to conjecture that any maximizing permutation must also be symmetric. Nevertheless, there is no general principle in mathematics that the solution of optimization problems involving a symmetric objective function must be symmetric. In fact, counterexamples to this "principle" abound. It is only by applying the information in Lemmas 1.8 and 1.9 in quite a subtle way that we can ultimately prove the Symmetry Lemma.

The second general observation is that Lemma 1.8 works by comparing a proposed maximizing permutation, $\pi$, with a "neighboring" permutation $\sigma = \pi\cdot(a\ b)$, where $(a\ b)$ is the 2-cycle in the permutation group $S_N$ that switches row $a$ and row $b$. We could call this a "local" comparison. But one problem with local comparisons is that they can leave us trapped at a local maximum. Symmetry gives us a way to get unstuck by potentially switching many rows at once.

Here is the central idea of the proof. Assume that $\pi$ is not symmetric about the AMD. Let $\pi'$ be the permutation obtained by reflecting $\pi$'s permutation matrix through the AMD. (We'll define this more formally below.) The strategy of the proof will be to find four other permutations, $\sigma$, $\tau_1$, $\tau_2$, and $\tau_3$, with the property that

$$F(\pi) + F(\pi') + F(\sigma) = F(\tau_1) + F(\tau_2) + F(\tau_3)$$

by inspection. The permutations $\sigma$ and $\tau_3$ will be related by a "good" row switch, which means we can apply Lemma 1.8 to conclude that $F(\sigma) > F(\tau_3)$. Hence

$$F(\pi) + F(\pi') < F(\tau_1) + F(\tau_2).$$

The remaining two pairs, $\pi$ and $\tau_1$ and $\pi'$ and $\tau_2$, are related by a complicated permutation (generally involving many row switches), so we cannot apply Lemma 1.8. However, we do know that $F(\pi) = F(\pi')$ by symmetry. Thus

$$F(\tau_1) + F(\tau_2) > 2F(\pi),$$

and this implies that either $F(\tau_1)$ or $F(\tau_2)$ is greater than $F(\pi)$. In either case, $\pi$ does not maximize the objective function. Thus (taking the contrapositive) any maximizer must be symmetric.

Here is a simple numerical illustration of the argument, taken from the case $N = 7$. We start with an asymmetric permutation $\pi$, illustrated below in green. We wish to show that it does not maximize the objective function, and we will do so by constructing two candidates $\tau_1$ and $\tau_2$, one of which must have a higher value. Note that $F(\pi) = 1 + 2640 + 176 + 1016 + 3256 + 3256 + 3312 = 13657$, so this is the number we are trying to beat. Also note that $\pi$ satisfies the Shape

Lemma (it is decreasing below the main diagonal and increasing above), so we cannot improve with any "local" move satisfying the hypotheses of Lemma 1.8.

| 1716 | 792 | 330 | 120 | 36 | 8 | 1 |
|------|------|------|------|------|------|------|
| 2640 | 1716 | 960 | 456 | 176 | 50 | 8 |
| 3102 | 2472 | 1716 | 1016 | 491 | 176 | 36 |
| 3312 | 2976 | 2416 | 1716 | 1016 | 456 | 120 |
| 3396 | 3256 | 2941 | 2416 | 1716 | 960 | 330 |
| 3424 | 3382 | 3256 | 2976 | 2472 | 1716 | 792 |
| 3431 | 3424 | 3396 | 3312 | 3102 | 2640 | 1716 |

The permutation π is given by π(1) = 7, π(2) = 1, π(3) = 6, π(4) = 5, π(5) = 2, π(6) = 3, and π(7) = 4. Thus the elements in positions (1, 7), (2, 1), (3, 6), (4, 5), (5, 2), (6, 3), and (7, 4) are highlighted above.

The first step is to work out the reflected permutation, π′, shown below in blue.

| 1716 | 792 | 330 | 120 | 36 | 8 | 1 |
|------|------|------|------|------|------|------|
| 2640 | 1716 | 960 | 456 | 176 | 50 | 8 |
| 3102 | 2472 | 1716 | 1016 | 491 | 176 | 36 |
| 3312 | 2976 | 2416 | 1716 | 1016 | 456 | 120 |
| 3396 | 3256 | 2941 | 2416 | 1716 | 960 | 330 |
| 3424 | 3382 | 3256 | 2976 | 2472 | 1716 | 792 |
| 3431 | 3424 | 3396 | 3312 | 3102 | 2640 | 1716 |

We see that π′(1) = 1, π′(2) = 5, π′(3) = 4, π′(4) = 1, π′(5) = 2, π′(6) = 3, and π′(7) = 6. Observe that π and π′ differ only in rows 2, 3, 4, and 7 (counting from the top) and columns 2, 3, 4, and 7 (counting from the right). They agree in rows 1, 5, and 6, so those rows will have no effect on our argument. I'll call them "irrelevant rows."

The next step is to work out the cycle structure of the permutation $\pi^{-1}\pi'$. It is easy to verify that $\pi^{-1}\pi'(2) = 4$, $\pi^{-1}\pi'(4) = 2$, $\pi^{-1}\pi'(3) = 7$, and $\pi^{-1}\pi'(7) = 3$. Thus $\pi^{-1}\pi'$ is a product of 2-cycles, namely (2 4)(3 7). This tells us that we are in Case 1 of the proof of the Symmetry Lemma, where $\pi^{-1}\pi'$ consists of two cycles, where one cycle contains the top row number (2) of the relevant rows, and the other contains the bottom row number (7). (Case 2, for any readers who are looking ahead, will be the case where $\pi^{-1}\pi'$ is a single cycle.)

The next step is to compute the permutations σ, $\tau_1$, $\tau_2$, and $\tau_3$. They are illustrated in the figures below. (All of them have the same entries in the irrelevant rows 1, 5, and 6, so only the entries in rows 2, 3, 4, and 7 are highlighted.) The first figure shows σ highlighted in yellow, along with π in green and π′ in blue. Note that σ simply consists of the entries along the anti-main diagonal in the set of relevant rows and columns. Also note that F(π) = F(π′) = 13657 and F(σ) = 12201.

i

| 1716 | 792 | 330 | 120 | 36 | 8 | 1 |
|------|-----|-----|------|------|------|------|
| 2640 | 1716 | 960 | 456 | 176 | 50 | 8 |
| 3102 | 2472 | 1716 | 1016 | 491 | 176 | 36 |
| 3312 | 2976 | 2416 | 1716 | 1016 | 456 | 120 |
| 3396 | 3256 | 2941 | 2416 | 1716 | 960 | 330 |
| 3424 | 3382 | 3256 | 2976 | 2472 | 1716 | 792 |
| 3431 | 3424 | 3396 | 3312 | 3102 | 2640 | 1716 |

Next we look at $\tau_1$ (gray), $\tau_2$ (red), and $\tau_3$ (pink). Note that *the very same elements of the* **P** *matrix are highlighted* in this figure and the previous one. This means we can be certain that

$$F(\pi) + F(\pi') + F(\sigma) = F(\tau_1) + F(\tau_2) + F(\tau_3).$$

In fact, $F(\tau_1) = 13678$, $F(\tau_2) = 13825$, and $F(\tau_3) = 12012$.

| 1716 | 792 | 330 | 120 | 36 | 8 | 1 |
|------|-----|-----|------|------|------|------|
| 2640 | 1716 | 960 | 456 | 176 | 50 | 8 |
| 3102 | 2472 | 1716 | 1016 | 491 | 176 | 36 |
| 3312 | 2976 | 2416 | 1716 | 1016 | 456 | 120 |
| 3396 | 3256 | 2941 | 2416 | 1716 | 960 | 330 |
| 3424 | 3382 | 3256 | 2976 | 2472 | 1716 | 792 |
| 3431 | 3424 | 3396 | 3312 | 3102 | 2640 | 1716 |

Note also that $\sigma$ and $\tau_3$ differ only by a "local" move, and we can apply Lemma 1.8 to conclude that $F(\sigma) > F(\tau_3)$.

It then follows by the argument above that either $F(\tau_1)$ or $F(\tau_2)$ must be strictly greater than $F(\pi)$. In this example, it turned out that they both were.

Note that in this example, the "best" permutation was $\tau_2$, which does not satisfy the Shape Lemma. (It has a decreasing sequence above the main diagonal.) But that's OK. The proof of the Symmetry Lemma merely returns a permutation that is more symmetric than $\pi$ (in this case, completely symmetric) and has a higher value of the objective function. It may take a few more applications of Lemma 1.8 and even a few more applications of the symmetrizing algorithm to get to the actual maximum.

One more point is important to make before we embark on the proof of the Symmetry Lemma. Although I have illustrated the argument using the matrix $\mathbf{P}^7$, *nothing in the proof requires the numerical values of* **P**. The whole proof applies equally well to any other mixed Monge matrix that is symmetric about the AMD. The theorem will therefore be stated in that generality.

To the reader wondering, "Where did the permutations $\tau_1$, $\tau_2$, and $\tau_3$ come from?" the general answer will be explained later. However, in the specific case where $\pi^{-1}\pi'$ is a product of two 2-cycles, as it was here, Figure 2.1 shows how they are constructed.

**Notation.** For any permutation $\pi \in S_N$, we define the permutation matrix representing $\pi$, denoted $M = M_\pi$, as follows: $M_{ij} = 1$ if $\pi(i) = j$, and $M_{ij} = 0$ otherwise. We will say that an ordered pair $(i, j)$ is "in the permutation matrix" if and only if $M_{ij} = 1$. The objective function $F$ is defined as $F(\pi) = M_\pi \cdot \mathbf{P}$, where $\mathbf{P}$ is any matrix that satisfies the conclusions of Lemma 1.4, Lemma 1.8 and Theorem 1.9 (i.e., a mixed Monge matrix that is symmetric about the anti-main diagonal). Note that $\mathbf{P}$ is symmetric about its anti-main diagonal (AMD). For any pair of coordinates $(i, j)$, its reflection in the AMD is given by $(N+1-j, N+1-i)$. Later it will be convenient to represent the "order-reversing" permutation by $\rho$; that is, $\rho(i) = N+1-i$. Given any permutation $\pi$, its reflection $\pi'$ is the permutation such that

$$\pi'(i) = j \leftrightarrow \pi(\rho(j)) = \rho(i) \leftrightarrow j = \rho\pi^{-1}\rho(i).$$

In other words, $\pi' = \rho\pi^{-1}\rho$. We will say that $\pi$ is symmetric (which means that $M_\pi$ is symmetric about the AMD) if and only if $\pi^{-1} \pi' = 1$, or equivalently $(\rho\pi)^2 = 1$.

Note also that there is a pointwise version of symmetry, which allows us to talk about symmetry of pairs of points in isolation from the rest of the permutation. A point $(i, \pi(i))$ lies on the AMD if and only if $\pi(i) = \rho(i)$, or $i$ is a fixed point of $\rho\pi$. Likewise, two points $(i, \pi(i))$ and $(j, \pi(j))$ form a symmetric pair with respect to the AMD iff $j = \rho\pi(i) \ne i$ and $\pi(j) = \rho(i)$. In this case, $(\rho\pi)^2(i) = \rho\pi(j) = i$. Thus a point $(i, \pi(i))$ and its reflection are both in the permutation matrix iff $i$ is a fixed point of $(\rho\pi)^2$ but not of $\rho\pi$.

**Lemma 2.1 (Symmetry Lemma).** If $\mathbf{P}$ is a mixed Monge matrix that is symmetric about the anti-main diagonal and if $\pi$ maximizes the objective function $F(\pi) = M_\pi \cdot \mathbf{P}$, then $M_\pi$ is also symmetric about the anti-main diagonal.

**Proof.** In order to show that $\pi$ is symmetric, we want to show that the nonzero entries of its permutation matrix either lie on the AMD or in a symmetric pair. To formalize this, let

$$T = T_\pi = \{i \le \lfloor \tfrac{N}{2} \rfloor \mid \rho\pi(i) = i\}$$

$$S = S_\pi = \{i \le \lfloor \tfrac{N}{2} \rfloor \mid \rho\pi(i) \ne i, (\rho\pi)^2(i) = i\}.$$

The notation is not accidental: T corresponds to "thrown" cards and S corresponds to "straight" cards. We will prove by induction that $S \cup T = \{1, \ldots, \lfloor N/2 \rfloor\}$. ($\lfloor N/2 \rfloor$ represents the integer part of $N/2$.) Technically this only shows that the top $\lfloor N/2 \rfloor$ rows and right $\lfloor N/2 \rfloor$ columns are symmetric about the AMD, but getting the rest will be very easy.

To do the induction, we will show that if $i \in S \cup T$ for all $i < a$, then $a \in S \cup T$. Actually, we will state the inductive hypothesis more precisely in a moment, but let's start with the first step.

**Step 1.** $1 \in S \cup T$. In fact, by the Shape Theorem we know already that for a maximizer of $F$ we must have $\pi(1) = N$. Thus $\rho\pi(1) = 1$, and $1 \in T$.

The idea of the inductive step is illustrated in Figure 2.2. Before we find $\pi(a)$, the first $(a-1)$ rows and the right $(a-1)$ columns of the permutation matrix are already occupied either by points on the AMD or by symmetric pairs. These points define the vertices of a nested sequence of squares that are each symmetric about the AMD, such that the remaining nonzero elements of the permutation matrix $M_\pi$ must lie inside the previously defined squares. Note also that the last $k(a-1)$ rows and the left $k(a-1)$ columns are also already occupied. Here $k(j)$ represents the number of symmetric pairs already found after $j$ steps; more precisely, $k(j) = \#(S \cap \{1, …, j\})$.

With this notation we can state our inductive assumption more precisely.

$IA_j$: $\{1, …, j\} \subset S \cup T$ and

$$\pi(\{1, …, j\} \cup \{N + 1 - k(j), …, N\}) = \{1, …, k(j)\} \cup \{N + 1 - j, …, N\}.$$

This is just a translation into symbols of what was already said about the nested squares. Notice that we have already shown that $1 \in T$ and hence that $k(1) = 0$. Hence the second part of $IA_1$ holds vacuously, adopting the convention that the set $\{1, …, 0\}$ denotes the integers greater than or equal to 1 and less than or equal to 0, i.e. the empty set. Thus we can go on to the inductive step.

**Step 2.** If $a \leq \lfloor N/2 \rfloor$, then $IA_{a-1}$ implies $IA_a$.

To prove this, we ask ourselves: What can $\pi(a)$ be equal to? By the inductive assumption $IA_{a-1}$, we know that $k(a-1) + 1 \leq \pi(a) \leq N+1-a$, because the left $k(a-1)$ columns and the right $(a-1)$ columns are already occupied. In fact, the inductive assumption implies that $k(a-1) + 1 \leq \pi(x) \leq N+1-a$ for *all* rows $x$ in the range from $a$ to $N-k(a-1)$. So for the remainder of the proof of Step 2 we can concentrate on this $(N - a - k(a-1))$ by $(N - a - k(a-1))$ submatrix of $M_\pi$, which is symmetric about the AMD. (It is the innermost of the nested matrices in Figure 2.2.)

If $\pi(a) = N+1-a = \rho(a)$, then clearly $a \in T$. Then $IA_a$ holds true, with $k(a) = k(a-1)$.

Now suppose $\pi(a) \neq N+1-a$, i.e., $(a, \pi(a))$ does not lie on the AMD. Likewise, the unique nonzero element of $M_\pi$ that lies in the $(N+1-a)$-th column is not on the AMD. This element has coordinates $(\pi^{-1}(N+1-a), N+1-a)$. These are two distinct nonzero entries of $M_\pi$ that form a decreasing sequence. Hence, by the Shape Theorem, one of them must lie below the main diagonal. Note that the latter point can also be written as $(\rho\pi'(a), \rho(a))$ because $\pi' = \rho\pi^{-1}\rho$. By a further reflection (noting that reflecting about the AMD maps points above the MD to points above the MD and points below the MD to points below the MD), we can conclude that either $(a, \pi(a))$ is below the main diagonal or $(a, \pi'(a))$ is below the main diagonal. Thus we can assume without loss of generality that $(a, \pi(a))$ is below the main diagonal, by relabeling $\pi$ as $\pi'$ if necessary.

Now we claim that in fact $\pi(a) = k(a-1) + 1$. If that were not the case, there would have to be some other row $c > a$ such that $\pi(c) = \underline{k}(a-1) + 1$. Then $(c, \pi(c))$ and $(a, \pi(a))$ would lie on a line with positive slope and yet they would be below the main diagonal, which violates the Shape Theorem again.

Thus by contradiction we conclude that $\pi(a) = k(a-1) + 1$, which places $(a, \pi(a))$ in the northwest corner of the inner submatrix shown in Figure 2.2.

Now let's ask about the bottom row of that submatrix, which is row $N - k(a-1)$. We will call this row $d$ for convenience. If $\pi(d) = \rho(a)$, then $(d, \pi(d))$ lies in the bottom right corner of the submatrix in Figure 3. Then $(a, \pi(a))$ and $(d, \pi(d))$ form a symmetric pair. We could then conclude that $a \in S$, we would increment $k$ by 1 (i.e., $k(a) = k(a-1) + 1$), and we would conclude that $IA_a$ is true.

So we have now gotten to the $64,000 question: What if $\pi(d) \neq \rho(a)$? This case is illustrated in Figure 2.3. We aim to show that this assumption leads to a contradiction.

One thing we do know is that $\pi'(d) = \rho(a)$, because the reflection of any point $(a, \pi(a))$ in the permutation matrix of $\pi$ is a point in the permutation matrix of $\pi'$. Hence we can start filling some information into Figure 2.3. The point $(a, \pi(a))$ lies in the northwest corner, the point $(d, \pi(d))$ on the bottom row but not in the southeast corner, and the point $(d, \pi'(d))$ is in the southeast corner. The permutation matrix for $\pi$ must contain a nonzero entry *somewhere* in the last column, so let's say it's in row $b$.

Now remember that $(a, \pi(a))$ lies below the main diagonal and therefore its symmetric twin $(d, \pi'(d))$ also lies below the MD. Because any point to the left of a point below the MD is also below the MD, $(d, \pi(d))$ lies below the MD as well. If $(b, \pi(b))$ were also below or on the main diagonal, then by Lemma 1.8 we could increase the value of the objective function by switching rows $b$ and $d$ in $M_\pi$. This would contradict the hypothesis that $\pi$ is optimal.

Hence, as illustrated in Figure 2.3, the point $(b, \pi(b)) = (b, \rho(a))$ lies strictly above the main diagonal. And therefore its reflection $(a, \rho(b))$ through the AMD also lies above the main diagonal. This point must also be equal to $(a, \pi'(a))$. [Once again, the reflection of a point in the permutation matrix of $\pi$ is *always* a point in the permutation matrix of $\pi'$.] So $\pi'(a) = \rho(b)$.

Notice that Figure 2.3 is starting to look a bit like Figure 2.1! This is no accident. We're just about ready to unveil the mysterious permutations $\sigma$, $\tau_1$, $\tau_2$, and $\tau_3$ that were mentioned earlier.

First, we note that $\pi^{-1}\pi'(d) = \pi^{-1}(\rho(a)) = b \neq d$. The permutation $\pi^{-1}\pi'$ can be written as a product of disjoint cycles, one of which must contain $d$. This cycle has $r \geq 2$ elements, so it can be written as $(d, b, a_1, \ldots, a_{r-2})$. For future reference, let's stipulate right now that there is no $a_0$. We adopt this convention because in the cases $r = 2$ and $r = 3$ it will turn out to be important to ignore certain vacuous equations involving $a_0$. (This is related to our earlier convention that the set $\{1, \ldots, 0\}$ denotes the empty set.) However, there will be no harm in identifying $d$ as being $a_{r-1}$ and $b$ as $a_r$, and we will at times do so.

Because $\pi^{-1}\pi'$ has an $r$-cycle of the form $(d, b, a_1, \ldots, a_{r-2})$, we have the following equations, where the third, fifth, etc. serve as definitions of $x_1, \ldots, x_{r-1}$:

$\pi'(d) = \rho a$

$\pi^{-1}(\rho a) = b$ or $\pi(b) = \rho a$

$\pi'(b) = x_1$

$\pi^{-1}(x_1) = a_1$ or $\pi(a_1) = x_1$

…

$\pi'(a_{r-2}) = x_{r-1}$

$\pi^{-1}(x_{r-1}) = d$ or $\pi(d) = x_{r-1}$.

By reflection through the AMD, we also conclude:

$\pi(a) = \rho d$

$\pi'(a) = \rho b$

$\pi(\rho x_1) = \rho b$

$\pi'(\rho x_1) = \rho a_1$

…

$\pi(\rho x_{r-1}) = \rho a_{r-2}$

$\pi'(\rho x_{r-1}) = \rho d$.

[If $r = 2$ the second-to-last equation is vacuous and the last equation is redundant.]

Observe from the above that $\pi^{-1}\pi'$ has another $r$-cycle, starting with $a$: $(a, \rho x_1, …, \rho x_{r-1})$. This is either identical to the previously mentioned cycle or disjoint. The two cases have slightly different proofs, and we will do the case where they are disjoint first.

**Case 1:** $\{a, \rho x_1, …, \rho x_{r-1}\} \cap \{d, b, a_1, …, a_{r-2}\} = \emptyset$.

Now we define the permutations $\sigma$, $\tau_1$, $\tau_2$, and $\tau_3$ as shown in Figure 2.4. Let's make a few comments to improve the reader's understanding. First of all, we are tacitly ignoring any rows other than $a, \rho x_1, …, \rho x_{r-1}, d, b, a_1, …, a_{r-2}$ and any columns other than $\rho a, x_1, …, x_{r-1}, \rho d, \rho b, \rho a_1, …, \rho a_{r-2}$. Other rows and columns relate to other cycles of $\pi^{-1}\pi'$ and we are free to modify $\pi$ and $\pi'$ as in Figure 2.4 without affecting any of those other rows and columns. Outside the rows and columns just listed, $\sigma$, $\tau_1$, $\tau_2$, and $\tau_3$ are all equal to $\pi$. (Remember the "irrelevant rows" and columns in our numerical example.)

We can describe $\sigma$ very simply as $\rho$ restricted to the $2r$ rows and columns we have listed. Likewise, $\tau_3$ is easy to describe: It is $\sigma$ modified by switching row $a$ and row $b$. This switch reduces the objective function because we have ensured that three of the vertices in the corresponding square lie strictly above the main diagonal.

The remaining two comparison permutations, $\tau_1$ and $\tau_2$, are harder to describe in words. I like to think of them as "recombined" versions of $\pi$ and $\pi'$, in which half of the permutation $\pi$ has been spliced together with half of the permutation $\pi'$ to form a new permutation. Furthermore, to continue the biological analogy, $\tau_1$ has been "transfected" with two bases from $\sigma$, namely the ones in row $a$ and row $b$.

Careful study of Figure 2.4 will confirm the following facts:

i)   The permutations $\tau_1$, $\tau_2$, and $\tau_3$ are in fact valid permutations, meaning that each row number and each column number occurs once and only once.

ii)  In each row $i$, the column elements $\pi(i)$, $\pi'(i)$, and $\sigma(i)$ are identical to the column elements $\tau_1(i)$, $\tau_2(i)$, and $\tau_3(i)$ in some order. This ensures that
$$M_\pi + M_{\pi'} + M_\sigma = M_{\tau_1} + M_{\tau_2} + M_{\tau_3}$$ and hence that
$$F(\pi) + F(\pi') + F(\sigma) = F(\tau_1) + F(\tau_2) + F(\tau_3).$$

iii) The permutations $\sigma$ and $\tau_3$ differ only in rows $a$ and $b$, so by Lemma 1.8, $F(\sigma) > F(\tau_3)$.

Just as in the numerical example, we then conclude that either $F(\tau_1) > F(\pi) = F(\pi')$ or that $F(\tau_2) > F(\pi) = F(\pi')$. This contradicts the assumption that $\pi$ and $\pi'$ are maximal for $F$.

**Case 2:** $\{a, \rho x_1, \ldots, \rho x_{r-1}\} = \{d, b, a_1, \ldots, a_{r-2}\}$.

In this case we have to use a slightly different version of the "recombine and transfect" argument. First, a few preliminaries. If the above two cycles are equal, then there must be some $m$ between 1 and $r-1$, (using the convention that $d = a_{r-1}$) such that $\rho x_1 = a_m$. Then, following the cycle around, we have $a_{m-1} = a$, $a_{m+1} = \rho x_2$, and so on up to $a_{r-2} = \rho x_{r-m-1}$, $d = \rho x_{r-m}$, $b = \rho x_{r-m+1}$, $a_1 = \rho x_{r-m+2}$, etc.

Also note that $\pi'(d) = \rho \pi^{-1} \rho(d) = \rho \pi^{-1}(x_{r-m}) = \rho(a_{r-m}) = x_{2r-2m+1}$ and similarly $\pi'(b) = x_{2r-2m+2}$, etc. (In all the above equations, subscripts should be interpreted modulo $r$.) Therefore $x_1 = x_{2r-2m+2}$. Hence $2r-2m+2 \equiv 1 \pmod{r}$. Because $1 \le m \le r-1$, the only possibility is that $r = 2m-1$ or, equivalently, $r$ is odd and $m = (r+1)/2$. This shows that cycles of $\pi^{-1}\pi'$ that are symmetric under the reflection $\rho$ cannot have a length that is an even number. But more importantly, the fact that $m = (r+1)/2$ gives us a way of "breaking the genome in half" and splicing the halves to form $\tau_1$ and $\tau_2$. Figure 2.5 shows how this is done.

The rest of the argument follows exactly as Case 1 did; observations (i), (ii) and (iii) hold and therefore we arrive at a contradiction of the assumption that $\pi$ and $\pi'$ are maximal for $F$.

Thus the answer to our "$64,000 question," posed several pages ago, is that $\pi(d)$ cannot be different from $\rho(a)$. Thus the inductive step is complete.

Finally, to complete the proof of the Symmetry Lemma, we note that we have shown that the first $\lfloor N/2 \rfloor$ rows and the right $\lfloor N/2 \rfloor$ columns of $\pi$ are symmetric about the AMD. In addition, some but not necessarily all of the bottom rows and left columns are also symmetric. This may leave us with a central square (as in Figure 2.2), bisected by the AMD, where the shape of $\pi$ is still undecided. The key point is that this central square lies entirely on or below the main diagonal, and hence the nonzero entries of $\pi$ within that square are strictly descending by the Shape Theorem. There is only one way to define a permutation matrix whose entries are strictly descending: it has to be the identity matrix. The identity matrix is, of course, symmetric about the AMD, and that completes the proof of symmetry. ∎

**A Few Final Remarks.**

(1) In the very last part of the argument, it is worth checking the case $r = 3$ (and $m = 2$) carefully, because this is the case where some vacuous equations occur in Figure 2.5. Figure 2.6 shows how the proof works (and the vacuous equations disappear) in that case.

(2) Note that in the case $r = 3$, there are exactly six possible maps from $\{b, a_1, a_2\}$ to $\{\rho a, x_1, x_2\}$. This fact explains, or at least hints at, why we need to consider six different permutations, $\pi$, $\pi'$, $\sigma$, $\tau_1$, $\tau_2$, and $\tau_3$. The most surprising thing to me was that these six permutations sufficed, not only for the case $r = 3$ but also for all larger values of $r$.

(3) Did we really need a proof by induction? The spirit of the proof is simpler than it appears. Basically we work out the cycle structure of $\pi^{-1}\pi'$ and then remove one cycle at a time (for cycles that are symmetric under $\rho$) or two at a time (for cycles that are not symmetric under $\rho$). The "removal" is accomplished by the splicing procedure illustrated in Figures 2.4 (for cycles that are not symmetric under $\rho$) and Figure 2.5 (for cycles that are symmetric under $\rho$). It appears to me that the cycles need to be removed in the correct order (from the "outside in," as shown in Figure 2.2) to guarantee that the objective function increases at each step. However, I am not entirely certain whether this is a real phenomenon or just an artifact of the proof.

Now let's step back and discuss what the Symmetry Lemma means in the context of the original Sun Bin's problem. Because of symmetry, each point $(i, \pi(i))$ in the optimal permutation matrix either lies on the AMD or else its symmetric partner $(\rho\pi(i), \rho(i))$ also lies in the permutation matrix. However, the line through any point not on the AMD and its symmetric partner is decreasing (in fact, parallel to the main diagonal). By the Shape Theorem, such symmetric pairs cannot occur above the main diagonal. Thus, every entry above the main diagonal in the optimal permutation matrix must lie on the AMD.

The points above the main diagonal correspond to cards in Stern-Mackenzie One Round War that are in set $T$ (they are "thrown" against higher-ranked cards of the opponent). This means that if $i \in T$, then $(i, \pi(i))$ lies on the AMD, and hence $\pi(i) = N + 1 - i$. In other words, if we decide to throw our $i$-th lowest card, then we must play it against the opponent's $i$-th highest card. As for the straight cards, they are simply played in order against all the opponent's remaining cards, so there is no choice about how to play them. Hence we have the following corollary:

**Corollary 2.2.** In Sun Bin's problem, an optimal permutation $\pi$ is uniquely determined by the set $T$ of thrown cards. Because $T \subset \{1, 2, \ldots, \lfloor N/2 \rfloor\}$, there are at most $2^{\lfloor N/2 \rfloor}$ remaining candidates for the optimal permutation.

## III. The No-Gaps Theorem

In view of the Symmetry Lemma, we now know that the maximizing permutation $\pi = \pi_{max}$ must have a very particular form. All the nonzero entries of $M_\pi$ above the main diagonal lie on the anti-main diagonal (AMD), and the remaining nonzero entries lie strictly below the main diagonal in a decreasing pattern that is symmetrical about the AMD. In fact, there is only one way to fill in the below-diagonal elements once the above-diagonal elements are known. Hence there are $2^{\lfloor N/2 \rfloor}$ candidates, each one indexed by a string of 1's and 0's, which are the elements of $M_\pi$ along the AMD, reading down from northeast to southwest. The meaning of the 1's and 0's is quite simple. A "1" in the i-th position of the string means that for the corresponding permutation, $\pi(i) = N + 1 - i$. A "0" means that $\pi(i) \neq N + 1 - i$.

For example, if $N = 9$, then the permutation matrix

$$\begin{bmatrix} 0 & 0 & 0 & 0 & 0 & 0 & 0 & 0 & 1 \\ 0 & 0 & 0 & 0 & 0 & 0 & 0 & 1 & 0 \\ 1 & 0 & 0 & 0 & 0 & 0 & 0 & 0 & 0 \\ 0 & 0 & 0 & 0 & 0 & 1 & 0 & 0 & 0 \\ 0 & 1 & 0 & 0 & 0 & 0 & 0 & 0 & 0 \\ 0 & 0 & 1 & 0 & 0 & 0 & 0 & 0 & 0 \\ 0 & 0 & 0 & 1 & 0 & 0 & 0 & 0 & 0 \\ 0 & 0 & 0 & 0 & 1 & 0 & 0 & 0 & 0 \\ 0 & 0 & 0 & 0 & 0 & 0 & 1 & 0 & 0 \end{bmatrix}$$

corresponds to the string (1, 1, 0, 1), and we can denote the corresponding permutation by $\pi_{1101}$. In the context of the original problem, this corresponds to the case where we throw our 1st, 2nd, and 4th weakest cards against the opponent's 1st, 2nd, and 4th strongest.

Intuitively, it is hard to imagine what advantage could be gained by throwing the 4th match instead of the 3rd. In other words, there should be no gaps between thrown tricks. This is the substance of the No Gaps Theorem. In reality, things are not so simple. For example, a "1" that is close to the main diagonal represents a card that is being thrown against an opponent's card that is not very much higher (e.g., our 4th weakest against his 5th weakest, above). This means we aren't throwing a whole trick; in the expectation, we are throwing only a little more than half a trick. Conceivably, under some circumstances, it might be a good strategy to give up half a trick to help out some of our "straight" cards (5, 6, 7, and 8 above)—and this might be better than giving up a whole trick to help 4, 5, 6, 7, 8, and 9. For sufficiently large $N$ this turns out not to work, but at least it's a plausible possibility that needs to be ruled out.

How large is "sufficiently large"? My most careful proof seems to be valid for $N > 400$ or so. I see little hope that it can be extended below about $N = 200$. However, it does not seem to be out of the realm of possibility that a computer could be programmed to verify the cases from $N = 1$ to 400, thus giving us a complete proof of the No Gaps Lemma (and therefore, also a complete solution to the problem).

In this section I will not attempt to give "best possible" estimates because this makes the logic of the proof harder to follow. I will simply give the simplest possible proof that is clearly valid for large enough $N$.

**Theorem 3.5. No Gaps Theorem.** For sufficiently large $N$, the optimal permutation corresponds to a string of 1's and 0's of the form (1, 1, …, 1, 0, …., 0), i.e., $k$ "1"s followed by $\lfloor N/2 \rfloor - k$ "0"s. Symbolically, $\pi_{max} = \pi_{11...10...0}$.

**Comment 1:** The No Gaps Theorem reduces the size of the set of possible solutions from $2^{\lfloor N/2 \rfloor}$ to $\lfloor N/2 \rfloor$. From the algorithmic point of view, this is arguably the most important step in the solution to Sun Bin's problem, because it reduces the complexity of a brute-force search from exponential time to linear time.

**Comment 2:** Note that the results of Sections 1 and 2 depended only on very general properties of the **P** matrix – the mixed Monge property and symmetry. That will not be the case in this section. We will have to do some very technical calculations involving binomial coefficients.

**Sketch of Proof:** There are two main ideas in the proof of the No Gaps Theorem: a local comparison and a global comparison. As discussed in Section 1, a local comparison is one that modifies the permutation matrix $M_\pi$ in a small number of rows and columns. In this section, we will modify three rows and columns.

Specifically, suppose $\pi$ is a permutation with a gap, such as $\pi_{1101}$ above. Let's call any 1's that appear after the first gap "rogue elements." The first of the main ideas in the proof is to compare this to the two permutations that are generated by sliding the first rogue element up one space or by sliding the last rogue element down one space. In other words, we will compare $\pi_{1101}$ to $\pi_{1110}$ and to $\pi_{1100}$. (In this second case we have committed an abuse of notation. In the case where the last rogue element is in the $\lfloor N/2 \rfloor$-th position, "sliding it down one place" means simply replacing it by a 0. Think of sliding a bead off the end of a string.)

The key claim is that one of these two comparison permutations is *guaranteed* to win more tricks than our candidate permutation $\pi$. So we can define an algorithm to find better and better permutations simply by repeatedly sliding the first rogue element up or the last rogue element down (and sometimes right off the end of the string). It is clear that this algorithm can only terminate when we get to a permutation with no rogue elements, i.e. a permutation that satisfies the No Gaps Theorem.

As noted, this is only one of the two main ideas. Now let me explain what goes wrong with the above argument, and why we also need to do a global comparison.

If we let $\sigma_1$ and $\sigma_2$ denote the two comparison permutations we have just defined, one by sliding the first rogue element of $\pi$ up and the other defined by sliding the last rogue element down, the matrices $M_{\sigma_1} - M_\pi$ and $M_\pi - M_{\sigma_2}$ have a specific form. They have only six nonzero elements, which lie at the vertices of an upside-down "L". Specifically, for some $i$ and $j$ (with $i > j$), they have a $-1$ in positions $(i, j)$, $(i+1, N-i)$, and $(N-j+1, N-i+1)$, and they have a $+1$ in positions $(i, N-

$i+1$), $(i+1, j)$, and $(N-j+1, N-i)$. Note that two vertices of the "L" are at positions $(i, N-i+1)$ and $(i+1, N-i)$, which lie on the AMD. (See Figure 3.1.)

Let's denote the matrix with the six nonzero entries in the six positions described by $N_{ij}$. Define $C_{ij} = N_{ij} \cdot \mathbf{P}$. Remember that $\mathbf{P}$ is the matrix that defines our objective function. In Sections 1 and 2 we used only general properties of $\mathbf{P}$, but in this section we will need the specific numerical values.

A particularly simple case, and a good one to use for motivation, is the case where there is only one rogue element, in row $(i+1)$. (See Figure 3.2.) In that case, $(M_{\sigma_1} - M_\pi) \cdot \mathbf{P} = C_{ij}$ and $(M_{\sigma_2} - M_\pi) \cdot \mathbf{P} = -C_{i+1,j+1}$. Our claim is that one of these numbers must be positive, which means that we can win more tricks by using permutation $\sigma_1$ or $\sigma_2$ rather than permutation $\pi$. A very simple way to show this would be to show that the $C_{ij}$'s are decreasing along diagonals parallel to the main diagonal, i.e.:

$$C_{ij} > C_{i+1,j+1}$$

for all $i, j$, because that would imply that $C_{ij} + (-C_{i+1,j+1}) > 0$.

Unfortunately, this is too optimistic. The statement in this generality is simply not true, and this is why we need the global comparison argument.

In the global comparison, we will compare $\pi$ to a very special "no-gaps" permutation, $\pi_0$, defined as follows:

$$\pi_0 = \pi_{11\ldots10\ldots0}$$

where the number of 1's in the subscript is exactly $\lfloor \sqrt{N \ln N / 2} \rfloor$.

At this point it looks as if the permutation $\pi_0$ has come completely from left field. However, there is a method behind our madness. As will become apparent both in this section and Section IV, **$\pi_0$ is, "up to first order," what the optimal permutation should be.** Curiously, except for $N = 2, 3, 4, 7$, and $8$, $\pi_0$ is *not* the optimal permutation, as we will show in Section IV. There are "higher order" effects that cause $\lfloor \sqrt{N \ln N / 2} \rfloor$ to be a slight overestimate of the number of tricks we should throw. However, as $N \to \infty$, this is *asymptotically* the correct number of matches to throw. So even though $\pi_0$ is not literally the optimal permutation, it is a very difficult target to beat, and we will learn a lot from this global comparison.

In particular, it will turn out that to beat $\pi_0$, we must throw *at most* $C\sqrt{N \ln N}$ tricks, for some $C$. Also, if there is a gap in the optimal permutation, the gap must be of length *at least* $N/4$. These statements translate into *a priori* inequalities, $i > N/4$ and $i - j < C\sqrt{N \ln N}$. Finally, we will show that *provided* $i > N/4$ and $i - j < C\sqrt{N \ln N}$ and *provided* that $N$ is "sufficiently large," the inequality $C_{ij} > C_{i+1,j+1}$ (or actually, a closely related inequality that is good enough for our purposes) is valid. In other words, it's not true for all $i$ and $j$, but it's true for the only values of $i$ and $j$ that matter.

Figure 3.3 illustrates the proof in the form of a flow chart. The input is a permutation $\pi = \pi_S$. Assuming that $S$ has a gap, the algorithm will do one of two things. It will either eliminate the first rogue element while strictly increasing the objective function, or it will produce a different permutation that strictly increases the objective function.

This completes the sketch of the proof of the No Gaps Theorem; now, let's get down to the details. We will begin with the part of the argument that deals with the special permutation $\pi_0$.

**Lemma 3.1.** In any diagonal of **P** that lies above the main diagonal, the maximum element(s) is (are) in the center. For any diagonal of **P** that lies below the main diagonal, the minimum element(s) is (are) in the center.

**Proof.** Keep in mind that $p_{i,i+k}$ counts the number of partitions $(\{a_1,...,a_N\},\{b_1,...,b_N\})$ of $\{1, 2, ..., 2N\}$ such that $a_i > b_{i+k}$. Let $U_{i,i+k}$ denote the set of such partitions. Furthermore, let
$U_{i,i+k}^1 \equiv$ {partitions such that $a_i > b_{i+k+1}$}
$U_{i,i+k}^2 \equiv$ {partitions such that $a_i + 1 = b_{i+k+1}$}
$U_{i,i+k}^3 \equiv$ {partitions such that $b_{i+k+1} > a_i + 1$}
$U_{i+1,i+k+1}^4 \equiv$ {partitions such that $b_{i+k} = 2i + k$ and $b_{i+k+1} = 2i + k + 1$}.
We leave it to the reader to show that $U_{i,i+k}$ is the disjoint union of $U_{i,i+k}^1$, $U_{i,i+k}^2$, and $U_{i,i+k}^3$, while $U_{i+1,i+k+1}$ is the disjoint union of $U_{i,i+k}^1$, $U_{i,i+k}^2$, and $U_{i+1,i+k+1}^4$. Hence, by a direct calculation, again left to the reader:

$$p_{i+1,i+k+1} - p_{i,i+k} = \#(U_{i+1,i+k+1}^4) - \#(U_{i,i+k}^3) = \frac{(2i+k-1)!(2N-2i-k-1)!}{i!(i+k)!(N-i)!(N-i-k)!} k(N - 2i - k).$$

By inspection, if $k > 0$ this last expression is positive when $N - 2i - k \geq 1$, or $i \leq (N - k - 1)/2$. If $k < 0$ the reverse is true. Note that if $i = (N - k)/2$, which can happen if $N-k$ is even, then $p_{i+1,i+k+1} - p_{i,i+k} = 0$, and this means that there are two maximum elements on the $k$-th diagonal, each one adjacent to the AMD. □

Now we turn to the "special" permutation $\pi_0 = \{N, N\text{-}1, \ldots, N\text{-}k+1, 1, 2, \ldots N\text{-}k\}$, where $k = \left\lfloor \sqrt{N \ln N / 2} \right\rfloor$. The next lemma gives a lower bound on the objective function applied to this permutation.

**Lemma 3.2.** If $N \geq 30$, then $M_{\pi_0} \cdot \mathbf{P} \geq (N - k)\binom{2N}{N}\left[1 - e^{-(k-1)^2/N}\right]$.

**Proof.** First, we can simply ignore the thrown tricks. That is, we use the lower bound $p_{1N} + p_{2,N-1} + ... + p_{k,N-k+1} \geq 0$. As for the non-thrown tricks, we note that they all lie on the $k$-th diagonal below the main diagonal, and the minimal value of $p_{ij}$ on this diagonal is in the center, per Lemma 1. Therefore
$$p_{k+1,1} + p_{k+2,2} + ... + p_{N,N-k} \geq (N - k)p_{i,i-k}$$
where $i = \lfloor (N + 1 - k)/2 \rfloor$.

We first consider the case where $(N-k+1)$ is even.

In this case, the minimum value $p_{\min} = p_{i,i-k}$ lies on the AMD. From Lemmas 1.2 and 1.3,
$p_{\min} = \binom{2N}{N} - \binom{N}{0}^2 - \binom{N}{1}^2 - \ldots - \binom{N}{j-1}^2$, where $j = (N-k+1)/2$.
Therefore we require a good upper bound for $\binom{N}{0}^2 + \binom{N}{1}^2 + \ldots + \binom{N}{j-1}^2$.

Lemma 3.8.2 in [LPV] gives an excellent upper bound for the similar sum without the squares. Repeating their proof essentially word for word with the squares gives

$$\binom{N}{0}^2 + \binom{N}{1}^2 + \ldots + \binom{N}{j-1}^2 \leq \frac{\binom{N}{j}^2}{2\binom{N}{\lfloor N/2 \rfloor}^2}\binom{2N}{N}.$$

At this point, we will diverge slightly from Lovasz et. al. because we have a better estimate for the numerator on the right-hand side. If $N$ is even, $N \geq 16$ and $3 \leq \frac{k-1}{2} \leq \frac{N}{4} - 1$, then Lemma 4.5 (which is independent of anything in this section) shows that
$\binom{N}{j} = \binom{N}{(N-k+1)/2} < \binom{N}{N/2}\exp[-(k-1)^2/2N]$. Plugging in $k = \lfloor\sqrt{N\ln N/2}\rfloor$, we find that that $(k-1)/2 > 3$ provided that $N \geq 30$. Hence, assuming $N \geq 30$, we can say that

$$\binom{N}{0}^2 + \binom{N}{1}^2 + \ldots + \binom{N}{j-1}^2 \leq \frac{1}{2}e^{-(k-1)^2/N}\binom{2N}{N}.$$

We will skip the details, but a similar statement holds if $N$ is odd, except that the constant "1/2" in the above inequality has to be replaced by $e^{1/N}/2$. This makes no practical difference because we will be interested in what happens for large values of $N$. For instance, because we have already assumed $N \geq 30$, we can replace "1/2" by $e^{1/30}/2$ or, even more simply, by 1. Thus

$$M_{\pi_0} \cdot \mathbf{P} \geq p_{k+1,1} + p_{k+2,2} + \ldots + p_{N,N-k} \geq (N-k)\binom{2N}{N}\left[1 - e^{-(k-1)^2/N}\right].$$

The same result holds if $N-k+1$ is odd, although the proof differs in minor details. Instead of one value of $p_{\min}$ lying on the AMD, there are now two adjacent minima lying just off the AMD. Instead of an exact formula for $p_{\min} = p_{i,i-k}$ we now have to use an inequality. From Lemma 1.8, $2p_{i,i-k} = p_{i,i-k} + p_{i+1,i-k+1} > p_{i,i+1-k} + p_{i+1,i-k}$. Each of the latter two terms lie on the AMD, and therefore we can represent them as a sum of squares of binomial coefficients. Finally, by a suitable modification of the Lovasz *et. al.* argument, we get the desired estimate on $M_{\pi_0} \cdot \mathbf{P}$. □

Lemma 2 has a very common-sense interpretation in terms of the percentage of tricks won by the second player in One Round War. Namely,

(proportion of tricks won) ≥
(proportion of tricks not thrown) × (lower bound for winning probability in tricks not thrown).

The proportion of tricks not thrown is, by definition, $1 - k/N$, and the lower bound for the winning probability in those tricks is $1 - e^{-(k-1)^2/N}$. Because $k \approx \sqrt{N\ln N/2}$, the proportion of tricks not thrown is roughly $1 - \sqrt{\ln N/2N}$. Also, $(k-1)^2/N \geq \frac{\ln N}{2}(1 - 4/\sqrt{N})$, from which it follows that

$$e^{-(k-1)^2/N} \leq \frac{1}{\sqrt{N}}\left(\sqrt{N}^{4/\sqrt{N}}\right) \leq e^{4/e}/\sqrt{N}.$$

Thus the lower bound for the winning probability in these tricks is at least $1 - e^{4/e}/\sqrt{N}$. The key thing to notice, disregarding the constants and the factor of $\ln N$, is that the proportion of tricks not thrown and the lower bound for the winning probability in those tricks are roughly equal. This explains the "mystery" of why the correct number of tricks to sacrifice is roughly $\sqrt{N \ln N / 2}$. If we sacrifice fewer tricks, say $\sqrt{N \ln N / c}$ for $c > 2$, then the probability of winning the non-thrown tricks is asymptotically $1 - C'N^{1/c}$, which is unacceptably low. If we throw too many tricks, then the proportion of tricks not thrown becomes unacceptably low. The constant $c = 2$ represents a happy compromise between throwing too many tricks and throwing too few. In Part IV we will confirm this heuristic argument that $\pi_0$ is near-optimal but not quite optimal.

Yet another way to interpret the above lemma is to estimate the proportion of tricks *lost*. If we let $p_L(\pi)$ denote the expected proportion of tricks lost when playing permutation $\pi$, then $p_L(\pi_0) < \sqrt{N}\left(e^{4/e} + \sqrt{\ln N / 2}\right)$. It is clear that, as $N \to \infty$, the dominant term in parentheses is the second one, and hence for sufficiently large $N$ we can say that $p_L(\pi_0) < \sqrt{N \ln N}$. It is also clear that, if we wanted, we could easily estimate how large $N$ has to be to make this statement true. However, at this point we are going to stop keeping track of what "sufficiently large" means, in order to make the argument easier to follow.

**Lemma 3.3 (a)** For sufficiently large $N$, the number of thrown tricks in the optimal permutation is less than $2\sqrt{N \ln N}$.
**(b)** If the optimal permutation, $\pi_{\max}$, has a gap, then the first rogue element of $M_\pi$ comes after row $\lfloor N/4 \rfloor$.

**Proof.** Part (a) is extremely easy. For each trick we throw, the expected number of wins is less than ½ and the expected number of losses is greater than ½. That's why we said we were "throwing" those tricks in the first place. If we throw $2\sqrt{N \ln N}$ tricks or more, the expected number of losses from *those tricks alone* is at least $\sqrt{N \ln N}$. That exceeds the expected number of losses from *all* tricks if we use permutation $\pi_0$. Hence it cannot be optimal to throw $2\sqrt{N \ln N}$ tricks or more.

For part (b), we will show that if the first rogue element comes too early, the objective function $F(\pi) = M_\pi \cdot \mathbf{P}$ can always be increased by moving the first rogue element to the previous row. This is a simpler version of the argument we will use in the proof of the main theorem, when we will either move the first rogue element up or the last rogue element down.

So let $\pi$ be a permutation with a gap, where the first rogue element appears in row $(i+1)$. Let $\sigma$ be the permutation defined by moving it up to row $i$. As discussed in the "Sketch of Proof," for some $j < i$,

$$F(\sigma) - F(\pi) = 2(p_{i+1,j} - p_{ij}) + p_{i,N+1-i} - p_{i+1,N-i} \equiv C_{ij}.$$

By applying Lemma 1.8 repeatedly, we know that $C_{ij} > C_{i1}$, and the latter is very easy to compute directly:
$$C_{i1} = 2(p_{i+1,1} - p_{i1}) + p_{i,N+1-i} - p_{i+1,N-i} = \binom{2N-i-1}{N-1} - \binom{N}{i}^2.$$
I claim that for large enough $N$, $\binom{N}{\lfloor N/4 \rfloor}^2 < \dfrac{1}{2}\binom{\lceil 7N/4 \rceil}{N}$. (*)

It will then follow that for $i \leq \lfloor N/4 \rfloor$,
$$\binom{N}{i}^2 \leq \binom{N}{\lfloor N/4 \rfloor}^2 < \frac{1}{2}\binom{\lceil 7N/4 \rceil}{N} \leq \frac{1}{2}\binom{2N-i}{N} < \binom{2N-i-1}{N-1},$$ and hence that $C_{i1}$ is positive.

The easiest case of (*) occurs when $N$ is divisible by 4. Then
$$\frac{\binom{7N/4}{N}}{\binom{N}{N/4}^2} = \frac{(7N/4)!(N/4)!^2(3N/4)!}{(N!)^3}.$$

To estimate this, we need a version of Stirling's formula with upper and lower bounds. A simple and adequate one is $\sqrt{2\pi N}\left(\dfrac{N}{e}\right)^N < N! < e\sqrt{N}\left(\dfrac{N}{e}\right)^N$. From this it follows that
$$\frac{\binom{7N/4}{N}}{\binom{N}{N/4}^2} \geq \frac{\pi^2\sqrt{21}}{4e^3} N^{1/2} \left[\frac{1}{2}\left(\frac{7}{4}\right)^{7/4}\left(\frac{3}{4}\right)^{3/4}\right]^N \geq \frac{1}{2}\sqrt{N}(1.07)^N.$$

The latter expression clearly increases without bound as $N \to \infty$, hence for large enough $N$ it is greater than 2.

The (similar) arguments in case $N = 4k+1$, $4k+2$, or $4k+3$ are left to the reader. $\square$

**Remark.** In part (b) of the preceding lemma, it would be delightful if we could improve the estimate to $\lfloor N/2 \rfloor$, because it would follow that there were no rogue elements at all (i.e., the No Gaps Theorem would be true). But there is not much wiggle room in the above proof, because the constant 1.07 at the end is already perilously close to 1. If it had been less than 1, then the key ratio would have been a decreasing exponential rather than an increasing exponential, and the proof would be invalid. Rough calculations suggest that part (b) holds with $\lfloor N/4 \rfloor$ replaced by $\lfloor 0.26N \rfloor$, but this tiny improvement would not be worth the extra effort.

We have one more lemma to prove before taking on the No Gaps Theorem. Earlier I compared the strategy to sliding beads on a string: you can always improve your permutation either by sliding the first rogue bead up or the last rogue bead down. But this strategy runs into difficulties in the case where the last rogue bead occurs in the last possible row, because there is nowhere to slide it down to. Lemma 3.4 says that if you have a bead (rogue or otherwise!) in the last position, you can improve your permutation simply by applying Lemma 1.8.

**Lemma 3.4.** Any permutation $\pi = \pi_{1**\ldots*1}$ with a 1 in the $\lfloor N/2 \rfloor$-th position on the anti-main diagonal is not optimal.

The proof of Lemma 4 looks more complicated than it is. Therefore I'll describe it first. If the last element on the AMD above the main diagonal is 1, then it forces the nonzero elements of $M_\pi$ that lie below the main diagonal to thread a very narrow needle. (That is, they must stay extremely close to the main diagonal.) If they don't thread the needle, then it is easy to improve the objective function by a simple switch, using Lemma 1.8.

In the case where $N$ is even, the needle is simply too narrow to be threaded. In the case where $N$ is odd, there is only one permutation that successfully threads the needle. That is the permutation $\pi = \pi_{100\ldots01}$. However, in this special case we use an argument involving Lemma 1.8 distributed over four rows and columns, to show that the permutation $\sigma = \pi_{100\ldots00}$ is better.

Now here is the formal version of the proof.

**Proof.** Case 1: $N$ is even. In this case, saying that there is a 1 in the $\lfloor N/2 \rfloor$-th position on the anti-main diagonal of $M_\pi$ means that $\pi(N/2) = N/2+1$. If $\pi(N/2+1)$ were larger than $N/2+1$, then $(N/2+1, \pi(N/2+1))$ would lie above the main diagonal, a contradiction. If $\pi(N/2+1)$ is smaller than $N/2$, then we switch the $N/2$ and $(N/2+1)$-th rows of $M_\pi$ to obtain a permutation matrix with a higher objective function, by Lemma 1.8.

That leaves the possibility that $\pi(N/2+1) = N/2$. Now consider the nonzero elements in the first $(N/2-1)$ columns of $M_\pi$. They cannot all lie in the first $(N/2-1)$ rows, because at least one of them would then be forced to lie on or above the main diagonal. In that case, by the symmetry lemma, such an element would have to lie on the AMD. But it can't both lie on the AMD and be in the upper left quadrant of $M_\pi$.

Thus, there is some column $j \leq (N/2-1)$ whose nonzero element lies in row $i > (N/2-1)$. This nonzero element cannot lie in row $N/2$ or $(N/2+1)$, because those rows are already taken. Hence $\pi(i) = j$ for some $i > (N/2+1)$ and some $j \leq (N/2-1)$. It follows that switching rows $(N/2+1)$ and $i$ in $M_\pi$ will increase the objective function, by Lemma 1.8. (Drawing a picture will make it easier to follow this argument.)

Case 2: $N$ is odd. Here, saying that there is 1 in the $\lfloor N/2 \rfloor$-th position on the anti-main diagonal of $M_\pi$ means that $\pi\left(\dfrac{N-1}{2}\right) = \dfrac{N+3}{2}$. Now think about $\pi\left(\dfrac{N+5}{2}\right)$. If $\pi\left(\dfrac{N+5}{2}\right) \geq \dfrac{N+5}{2}$, we would have a descending pair of 1's on or above the main diagonal, in violation of Theorem 1.9. If $\pi\left(\dfrac{N+5}{2}\right) \leq \dfrac{N-1}{2}$, we could increase the objective function by switching rows $(N-1)/2$ and $(N+5)/2$. (Looking at Figure 3.4b and comparing to Lemma 1.8 will make this clear.) Since $\dfrac{N+3}{2}$ is taken, the only remaining possibility is $\pi\left(\dfrac{N+5}{2}\right) = \dfrac{N+1}{2}$. By similar arguments with Lemma 1.8, we can show that $\pi\left(\dfrac{N+3}{2}\right)$ must be $\dfrac{N-1}{2}$ and $\pi\left(\dfrac{N+1}{2}\right)$ must be $\dfrac{N-3}{2}$. The remaining entries in the first $(N-5)/2$ columns must fall in the first $(N-3)/2$ rows, by an argument like the one in Case 1. The only way to arrange this, and have them all below the main diagonal, is to have all these entries lie on the first sub-main diagonal. By the Symmetry Lemma, all the remaining entries are uniquely defined, and we find that $\pi = \pi_{10\ldots01}$. Note that this permutation differs from $\sigma = \pi_{10\ldots00}$ only in the four rows and columns pictured in Figure 3.4b. Hence, for the rest of the argument we can look only at this $4 \times 4$ submatrix.

To complete the argument, I claim that this permutation has a worse objective function than σ does. In other words, we can "slide the bead off the end of the string." The argument is so clumsy to write down symbolically that I will use a pictorial proof (Figure 3.4c). The first equality follows from Lemma 1.3; the second equality follows from the fact that all the entries on the main diagonal of **P** are equal; and the third inequality follows from Lemma 1.8. □

**Remark.** It's possible that "sliding the bead off the end of the string" always works, i.e. the permutation $\pi_{xxx...x1}$ is always worse than the permutation $\pi_{xxx...x0}$. However, we have not needed to prove such a general statement. We merely proved that $\pi_{xxx...x1}$ is worse than *something*, except in the specific case of the permutation $\pi_{100...01}$ when $N$ is odd. And in that special case, $\pi_{100...01}$ indeed turns out to be worse than $\pi_{100...00}$.

**Proof of No Gaps Theorem.** Recall that the theorem says that the optimal permutation has no gaps. Therefore, let us suppose that $\pi = \pi_{1...0...1...}$ is an eligible permutation (i.e., satifying the shape lemma and symmetry lemma) with a gap. We wish to prove that $\pi$ is not optimal.

Assume that the first rogue element lies in row $(i+1)$. That is, the string $S = (1, ..., 0, ..., 1, ...)$ starts out with $m$ 1's, followed by $(i-m)$ 0's, followed by $r$ 1's. We will define two comparison strings. $S_1$ starts with $m$ 1's followed by $(i-m-1)$ 0's, followed by a 1, a 0, and $(r-1)$ 1's. (This corresponds to "sliding the first bead up one row.") $S_2$ is the string defined by sliding the "1" in position $(i+r)$ down one row. Note that, by Lemma 4, there must be a place to slide it to, because the string $S$ cannot end with a "1." (If it did, then $\pi$ would not be optimal, and then we would be done.) Finally, let $\sigma_1 = \pi_{S_1}$ and $\sigma_2 = \pi_{S_2}$. These are our two comparison permutations.

As noted in the preamble, $\left(M_{\sigma_1} - M_\pi\right) \cdot \mathbf{P} = C_{ij}$, where $j = i - m$. We note that by Lemma 3, $i > \lfloor N/4 \rfloor$ and $i - j = m < 2\sqrt{N \ln N}$. Also, $\left(M_\pi - M_{\sigma_2}\right) \cdot \mathbf{P} = -C_{i+r, j+1}$. (To see this, notice that the column number $j$ does not increment when we go through a consecutive string of "rogue 1's." It only increments when we get to a "0" in the string $S$.)

I now claim that $C_{ij} \leq 0 \Rightarrow C_{i+1, j+1} < C_{ij} \leq 0$ (**). In the simple case where $r = 1$ (i.e., an isolated rogue element) this immediately means that if $\sigma_1$ is worse than $\pi$ (i.e., if $M_{\sigma_1} \cdot \mathbf{P} < M_\pi \cdot \mathbf{P}$) then $\sigma_2$ is better than $\pi$, and hence $\pi$ is not optimal. The case for general $r$ is not much harder, and we will do it at the end.

To prove the claim (**) we will have to get our hands dirty and do some very technical calculations with binomial coefficients. By definition and symmetry about the AMD,

$$C_{ij} = 2(p_{i+1,j} - p_{ij}) + (p_{N+1-i} - p_{i+1,N-i}) = 2\binom{2N-i-j}{N-i}\binom{i+j-1}{i} - \binom{N}{i}^2 = \frac{2j}{i+j}\binom{2N-i-j}{N-i}\binom{i+j}{i} - \binom{N}{i}^2.$$

We apply the following known identity involving binomial coefficients:

$$\binom{m}{k} / \binom{n}{k} = \binom{n-k}{m-k} / \binom{n}{m},$$

first with $m = i+j$, $n = N$, and $k = i$ to obtain

$$\binom{i+j}{i}\bigg/\binom{N}{i} = \binom{N-i}{j}\bigg/\binom{N}{i+j}$$

and secondly with $m = N$, $n = 2N - j$ and $k = i$ to obtain

$$\binom{N}{i}\bigg/\binom{2N-j}{i} = \binom{2N-i-j}{N-i}\bigg/\binom{2N-j}{N}.$$

Therefore

$$\binom{2N-i-j}{N-i}\binom{i+j}{i} = \binom{N}{i}^2\binom{N-i}{j}\binom{2N-j}{N}\bigg/\binom{N}{i+j}\binom{2N-j}{i} = \binom{N}{i}^2\binom{2N}{N}\binom{N}{j}\bigg/\binom{2N}{i+j}\binom{N}{i}.$$

So

$$C_{ij} = \left[2\binom{2N}{N}\frac{j}{i+j}\frac{\binom{N}{j}}{\binom{2N}{i+j}\binom{N}{i}} - 1\right]\binom{N}{i}^2 \equiv \left[2\binom{2N}{N}D_{ij} - 1\right]\binom{N}{i}^2,$$

where we are defining $D_{ij} = \dfrac{j\binom{N}{j}}{(i+j)\binom{2N}{i+j}\binom{N}{i}}$. I claim that $D_{ij} > D_{i+1,j+1}$ (***). By a direct and extremely tedious calculation,

$$D_{ij} - D_{i+1,j+1} = \frac{i!(N-i-1)!(i+j-1)!(2N-i-j-2)!}{(2N)!\,j!(N-j)!} \times$$

$$\times \left[2j(2N+1)(N-j)(N-2i-1) - (i-j)^2(i-j+1)N - (i-j)(N-j)\right]. \text{ (****)}$$

We merely need to show that the expression in braces (****) is positive. To convey the idea of the argument most simply, let's just look at the highest-order terms (the terms of degree four). These terms are $4jN(N-j)(N-2i) - (i-j)^3 N$. Now is when we will invoke the *a priori* inequalities involving $i$ and $j$. First,

$$j = i - (i-j) \geq \frac{N}{4} - 2\sqrt{N\ln N} > \frac{N}{8},$$

where the latter inequality holds for "sufficiently large $N$," and we will leave it to the reader to worry about how large that means. Also, of course, $j < N/2$ because $i < N/2$.

Thus the fourth-order terms are bounded below by

$$4\left(\frac{N}{8}\right)\left(\frac{N}{2}\right)(N - 2i) - 8N^{5/2}(\ln N)^{3/2} = \frac{N^3}{4}(N - 2i - 1) + \frac{N^3}{4} - 8N^{5/2}(\ln N)^{3/2} > 0,$$

where the last inequality holds for sufficiently large $N$. As an aside, we point out that it takes quite a long time for $8N^{5/2}(\ln N)^{3/2}$ to become smaller in absolute value than $N^3/4$. Taking $N > 10^7$ will work. This is the step that is most in need of improvement if we want to establish the No Gaps Theorem for smaller values of $N$.

Hence for sufficiently large $N$, $D_{ij} > D_{i+1,j+1}$ as claimed (***). The above argument would have failed if $N = 2i$, because the first two terms would have canceled. However, that is impossible because it would mean that the first rogue element is in row $N/2+1$, which is impossible because the string $S$ has only $\lfloor N/2 \rfloor$ terms.

Because we looked only at the highest-order terms in expression (****), this proof is not really complete. However, the lower-order terms can be handled in exactly the same way (the only change is that the criterion for "sufficiently large" becomes even larger). The details are left to the reader.

To show claim (**) observe that if $C_{ij} \leq 0$, then

$$C_{i+1,j+1} = [2\binom{2N}{N}D_{i+1,j+1} - 1]\binom{N}{i+1}^2 < [2\binom{2N}{N}D_{ij} - 1]\binom{N}{i+1}^2 \leq [2\binom{2N}{N}D_{ij} - 1]\binom{N}{i}^2 = C_{ij},$$

where the last inequality holds because $\binom{N}{i+1} \geq \binom{N}{i}$ and because $2\binom{2N}{N}D_{ij} - 1$ is nonpositive.

Finally, going back to the paragraph before claim (**) was made, we know that $C_{i+r,j+1} < C_{i+r,j+r}$ from Lemma 1.8, and hence

$$-C_{i+r,j+1} > -C_{i+r,j+r} > \ldots > -C_{i+1,j+1} > -C_{ij} \geq 0.$$

Therefore if $\sigma_1$ is not better than $\pi$, then $\sigma_2$ is better than $\pi$. Thus, $\pi$ is not optimal, and the proof of the No Gaps Theorem is complete. $\square$

**Remarks.** All steps in the above proof are valid for $N$ greater than 10 million. *Very* substantial improvements are possible.

In Lemma 2, the key idea was to give a lower bound for the sum of $p_{ij}$'s along the $k$-th diagonal by taking $(N-k)$ times the smallest element. Instead, in Part IV we will develop a precise formula for the sum of these $p_{ij}$'s, and by using this formula we can give a sharper lower bound on the number of tricks won in the "special" permutation $\pi_0$. This leads to a corresponding reduction in the estimate of the "threshold" for $N$ (i.e., when $N$ becomes sufficiently large for the No Gaps Theorem to be true.)

The proof of Lemma 3(a) can also be improved a great deal. Our lower bound on the number of tricks lost was based on the idea that we have a probability at least ½ of losing each of the "thrown" tricks. In fact, the probability of losing them is normally quite close to 1. Again, using a better estimate of these probabilities reduces the threshold for $N$.

In the proof of the No Gaps Theorem itself, the threshold for $N$ turns out to be quite sensitive to a lower bound on $N - 2i$. In the proof above we simply used the fact that $N - 2i \geq 1$. In fact, using Lemma 3.4 we could already have said that $N - 2i \geq 4$. If this inequality could be improved to $N - 2i \geq 10$, then (in conjunction with all the other improvements just mentioned) I could show that the No Gaps Lemma holds for $N > 400$.

## IV. The Number of Tricks to Throw

In the previous sections we have gradually narrowed down the form of the optimal permutation in Sun Bin's problem. In Part I, the "mixed Monge property" was used to prove that the optimal strategy for player 2 is to divide his cards into two groups. One group, $T$, consists of cards that will intentionally be paired against higher-ranked cards of the opponent. These are the tricks that we "throw." The second group, $S$, consists of cards that will intentionally be paired against lower-ranked opponents. These are the cards that are played "straight."

In Part III we showed that the optimal strategy has no gaps in it if $N$ is sufficiently large (in particular, if $N > 10^7$). Thus, for instance, it could not be optimal to choose $T = \{1, 2, 4\}$. If we decide to throw three tricks, then we should always throw tricks 1, 2, and 3. More generally, $T$ must have the form $\{1, 2, \ldots, k\}$ for some $k < N/2$. At this point we have reduced the number of possible candidates for the optimal permutation all the way from $N!$ to $\lfloor N/2 \rfloor$.

The objective of Part IV is to derive an exact description of the optimal number of tricks to throw, *assuming the No Gaps Conjecture is true*. In light of Part III, this means that the results in this section are unconditionally true for $N$ greater than 10 million. I conjecture that they are true for all $N$. Define $k^*(N)$ to be the number of elements of $T$ in the optimal permutation. The main result of this section is as follows:

**Theorem 4.1. (Number of Tricks Theorem)** For any $N$ for which the No Gaps Conjecture holds,

$$k^*(N) = \sup\left\{k : \binom{N}{k-1}^2 + \sum_{j=0}^{N-k}\binom{2N}{j} \geq \binom{2N}{N}\right\}.$$

In particular, this formula is valid for $N$ greater than 10 million.

It should be noted that the first term in the above expression, $\binom{N}{k-1}^2$, accounts for the expected number of tricks that we attempt to "throw" but actually win. This term is much smaller than the expected number of wins in the non-thrown tricks. The formula would be simpler, and the mathematical analysis likewise, if this "nuisance term" were not present. This would correspond to literally forfeiting the thrown tricks instead of playing them and hoping for the best. Thus we will do most of our analysis on the simpler function

$$k(N) = \sup\left\{k : \sum_{j=0}^{N-k}\binom{2N}{j} \geq \binom{2N}{N}\right\}$$

instead. This is the number of the sub-main diagonal of **P** with the largest sum; for example, if the diagonal one step below the main diagonal has the largest sum, then $k(N) = 1$. The difference between $k^*(N)$ and $k(N)$ can never be greater than 1, and I do not know a single value of $N$ where $k^*$ and $k$ are different. Thus, here is a simpler strategy that is optimal "for all practical purposes."

**"Practical Version" of Theorem 4.1.** For any $N$ for which the No Gaps Conjecture holds, you should pair your cards $\{1, 2, 3, \ldots, k(N)\}$ against the opponent's cards $\{N, N-1, \ldots, N+1-k(N)\}$. To compute $k(N)$ go to your school's mathematics department and ask for a printout of the $2N$-th row of Pascal's triangle. Then add up the numbers in the $2N$-th row, from left to right, until the sum is greater than the center element in that row. Then $k(N)$ is $N+1$ minus the number of elements you have just added up. (For two examples, see the Introduction.)

**Alternative Method.** We can also compute $k(N)$ by adding from the "center out":
$$k(N) = \sup\left\{k : \left(\frac{1}{2}\right)4^N - \left(\frac{3}{2}\right)\binom{2N}{N} - \binom{2N}{N-1} - \cdots - \binom{2N}{N-k+1} \geq 0\right\}.$$
While this formula is a little less elegant, it would be quicker for practical use because it involves adding only $k$ rather than $N$-$k$ numbers (and in practice, $k$ is on the order of $\sqrt{N}$, as discussed below).

**Explicit Estimates.** It is natural to wonder how rapidly the functions $k(N)$ or $k^*(N)$ grow as $N \to \infty$. We can prove that
$$\sqrt{N\ln(N)/4} < k(N) < \sqrt{N\ln(N)/2}$$
for all $N \geq 400$. Computer calculations (see Table 4.1) show that these inequalities also hold for $399 \geq N \geq 91$, so that $N = 90$ is the last exception. (When $N = 90$, the optimal number of tricks to throw is 10, which is slightly less than $\sqrt{90\ln(90)}/2 = 10.06...$)

In addition, we can prove that asymptotically,
$$k(N) \sim \sqrt{N\ln(N)/2}.$$
Interestingly, $k(N)$ appears to converge to this limit very slowly. Even at $N = 500$, computer calculations show that $k(N)$ is still considerably closer to $\sqrt{N\ln(N)/4}$ than to $\sqrt{N\ln(N)/2}$.

The same asymptotic estimate holds, of course, for $k^*(N)$ because $k(N) \leq k^*(N) \leq k(N)+1$.

We will now prove Theorem 1. The first step is to derive a recursive formula for $\mathbf{P}^N$.

**Theorem 4.2. (Hook-Sum Formula)** For any $N$, $I$, $K$, such that $1 \leq I, K \leq N$,
$$p_{I,K+1}^{N+1} = \sum_{j=K}^{N}(p_{Ij}^N - p_{I-1,j+1}^N) + \sum_{j=1}^{I-1}(p_{jK}^N - p_{j-1,K+1}^N).$$
Before moving on to the proof of this theorem, it will be very helpful to look at an example. (Also see Table 1.1, which gives two more examples.) Consider the fourth matrix,
$$\mathbf{P}^4 = \begin{bmatrix} 35 & 15 & 5 & 1 \\ 55 & 35 & 17 & 5 \\ 65 & 53 & 35 & 15 \\ 69 & 65 & 55 & 35 \end{bmatrix}.$$
Suppose we want to compute the (2, 3)-th entry of $\mathbf{P}^5$, which is 66. According to Theorem 2, we highlight the (2, 2)-th entry of $\mathbf{P}^4$, which is 35. We add this element, all the elements to its right, and all the elements above it. These elements form a right-angled figure that looks like a hook, so I call this sum a hook-sum. Now we *subtract* the hook-sum of the next entry to the northeast of 35, which is 5. The difference of these two hook-sums is the (2, 3)-th entry of $\mathbf{P}^5$:
$$p_{23}^5 = (35 + 17 + 5 + 15) - (5 + 1) = 72 - 6 = 66.$$
For elements in the first row or last column, the second hook-sum is zero (there is no entry northeast of the highlighted entry). Thus, for instance,
$$p_{13}^5 = 15 + 5 + 1 = 21.$$

In this way we can derive the 4-by-4 matrix that forms the northeast corner of $\mathbf{P}^5$:

$$\mathbf{P}^5 = \begin{bmatrix} 126 & 56 & 21 & 6 & 1 \\ 196 & 126 & 66 & 26 & 6 \\ 231 & 186 & 126 & 66 & 21 \\ 246 & 226 & 186 & 126 & 56 \\ 251 & 246 & 231 & 196 & 126 \end{bmatrix}.$$

To derive the remaining entries (the first column and last row) we could of course, simply use the symmetry properties of **P**. But a more self-contained, and I think more elegant, procedure is to "augment" $\mathbf{P}^N$ by adding a first column and last row whose entries are all equal to $\binom{2N}{N}$. We then perform the hook-sum procedure on this augmented matrix. Here is how the first few steps of the recursion work:

$\mathbf{P}^1 = [1]$. Augment: $\hat{P}^1 = \begin{bmatrix} 2 & 1 \\ 2 & 2 \end{bmatrix}$.

Hook-sum: $\mathbf{P}^2 = \begin{bmatrix} 3 & 1 \\ 5 & 3 \end{bmatrix}$. Augment: $\hat{P}^3 = \begin{bmatrix} 6 & 3 & 1 \\ 6 & 5 & 3 \\ 6 & 6 & 6 \end{bmatrix}$.

Hook-sum: $\mathbf{P}^3 = \begin{bmatrix} 10 & 4 & 1 \\ 16 & 10 & 4 \\ 19 & 16 & 10 \end{bmatrix}$. Augment: $\hat{P}^4 = \begin{bmatrix} 20 & 10 & 4 & 1 \\ 20 & 16 & 10 & 4 \\ 20 & 19 & 16 & 10 \\ 20 & 20 & 20 & 20 \end{bmatrix}$, etc.

This procedure can be repeated for as many iterations as desired. More formally, Let us define a basis $\{H_{ij}^N\}$ for the vector space of N-by-N matrices as follows: Let $H_{ij}^N$ be the matrix which has 1 in each position in the ij-th hook-sum, -1 in each position in the (i-1, j+1)-th hook-sum, and 0 everywhere else. For example,

$$H_{32}^4 = \begin{bmatrix} 0 & 1 & -1 & 0 \\ 0 & 1 & -1 & 0 \\ 0 & 1 & 1 & 1 \\ 0 & 0 & 0 & 0 \end{bmatrix}.$$

Define the "hook-sum map" of any N-by-N matrix, $\mathcal{H}(M)$, to be the N-by-N matrix $M'$ whose (ij)-th entry is $H_{ij}^N \cdot M$.

Define the augmented matrix $\hat{P}^N$ to be the (N+1) × (N+1) matrix whose northeast N × N submatrix is $\mathbf{P}^N$, and whose first column and last row consist entirely of the number $\binom{2N}{N}$ repeated. Then Theorem 2 can be restated as follows:

$$\mathbf{P}^{N+1} = \mathcal{H}(\hat{P}^N).$$

To be more precise, Theorem 2 shows that the northeast $N \times N$ submatrices of both sides are equal. Verifying that the first column and the last row agree requires a separate calculation, which is left to the reader. (Lemma 1.6 is an important ingredient in the proof.)

To prove Theorem 2, we will first require two small lemmas.

**Lemma 4.3 (Also mentioned in Section 1 as Lemma 1.2)**
$$p_{ij}^N = \sum_{k=1}^{i} \binom{i+j-1}{k-1}\binom{2N-i-j+1}{N-k+1}.$$
As discussed in Section 1, this formula is useful because it gives us a simple expression for the difference of consecutive terms along a SW-NE diagonal:
$$p_{ij} - p_{i-1,j+1} = \binom{i+j-1}{i-1}\binom{2N-i-j+1}{N-i+1}.$$
This is obviously relevant to Theorem 2 because the hook-sum map can be expressed as a sum of terms of this form. However, I emphasize that *both* Lemma 1.1 *and* Lemma 4.3 are essential to the proof of Theorem 4.2.

**Proof of Lemma 3.** We have a nice synthetic proof that doesn't involve any calculations. Recall that $p_{ij}$ counts the number of ways to partition $2N$ cards, $\{1, \ldots, 2N\}$ into two hands, $\{a_1, \ldots, a_N\}$ and $\{b_1, \ldots, b_N\}$ (with players ordered from lowest to highest) in such a way that $a_i > b_j$. Let $S_{ij}$ denote the set of partitions $(\{a_1, \ldots, a_N\}, \{b_1, \ldots, b_N\})$ such that $a_i > b_j$. Let $R_{ij}$ denote the set of partitions $(\{a_1, \ldots, a_N\}, \{b_1, \ldots, b_N\})$ such that $\{a_1, \ldots, a_i, b_1, \ldots b_j\} = \{1, \ldots, i+j\}$. In other words, the elements of $R_{ij}$ are partitions such that the $i$ weakest cards of hand 1 and the $j$ weakest cards of hand 2 collectively form the $(i+j)$ weakest cards of the combined hands.

It is then a simple exercise to prove that the set-theoretic difference $S_{ij} \setminus S_{i-1,j+1} = R_{i-1,j}$. The elements in the left-hand set are counted by $p_{ij} - p_{i-1,j+1}$. The elements in the right-hand set are counted by the product $\binom{i+j-1}{i-1}\binom{2N-i-j+1}{N-i+1}$, because we can generate all the partitions in $R_{ij}$ first by choosing the weakest $(i-1)$ players on team 1 from $\{1, \ldots, i+j-1\}$ and then by choosing the strongest $N-i+1$ players from $\{i+j, \ldots, 2N\}$. □

**Lemma 4.4.** For $1 \leq I \leq N$,
$$\sum_{j=1}^{I-1} (p_{j,K}^N - 2p_{j,K+1}^N + p_{j,K+2}^N) = \binom{I+K-1}{I-2}\binom{2N-I-K}{N-I}.$$
Unlike Lemma 4.3, I don't know a good combinatorial interpretation of this equation; it's just an identity we will need in the proof of Theorem 4.2.

**Proof of Lemma 4.4.** Note that the term in parentheses on the left-hand side is a second difference of the **P** matrix along the rows. Hence, by using Lemma 1.1, it is easily represented in terms of binomial coefficients. The formula can then be proved in routine but tedious fashion by induction on $I$, and the details are left to the reader. □

**Proof of Theorem 4.2.** This is again a proof by induction, this time on $K$. It is straightforward to show that the first row and last column of $\mathbf{P}^N$ satisfy the indicated formula. We then induct backward across the rows. Assume that $p_{I,K+2}^{N+1}$ satisfies the identity in the statement of the theorem (with $K$ replaced by $K+1$, of course). We then want to show that $p_{I,K+1}^{N+1}$ satisfies the same identity. We write
$$p_{I,K+1}^{N+1} = p_{I,K+2}^{N+1} + (p_{I,K+1}^{N+1} - p_{I,K+2}^{N+1}).$$
The first term on the right-hand side can be evaluated by the inductive assumption. The second term on the right-hand side is a difference of two adjacent elements of $\mathbf{P}^{N+1}$ along a row. Hence it is easily evaluated by applying the definition of $\mathbf{P}^{N+1}$.

Having done this, we then evaluate the difference
$$p_{I,K+1}^{N+1} - \sum_{j=K}^{N}(p_{Ij}^N - p_{I-1,j+1}^N) - \sum_{j=1}^{I-1}(p_{j,K}^N - p_{j-1,K+1}^N).$$
Of course, the goal is to prove that this difference equals 0. The main point is that the calculation involves a difference of two hook-sum maps; we end up computing $\mathbf{P}^N \cdot (H_{IK}^N - H_{I,K+1}^N)$.

In this difference of hook-sums, most of the entries along the $I$-th row cancel. The entries along the $K$-th, $(K+1)$-th, and $(K+2)$-th columns do *not* cancel, but form a series of "second differences" that is, in fact, precisely the left-hand side of the equation in Lemma 4. Finally, there are three leftover terms that arise around the corner of the $(I, K)$-th hook sum as well as from the inductive assumption. These leftover terms can be shown to add to exactly the right-hand side of the equation in Lemma 4.4. By Lemma 4.4, the difference between the left-hand side and the right-hand side is zero, and that completes the inductive proof. (The nitty-gritty details of the calculation are left to the reader.) □

Now we proceed to the proof of Theorem 4.1 (the Number of Tricks Theorem).

*Proof of Theorem 4.1.* Let $s_{kN}$ denote the sum of the $k$-th diagonal above the main diagonal in $\mathbf{P}^N$, and let $\bar{s}_{kN}$ denote the sum of the $k$-th diagonal below the main diagonal. Thus, $s_{kN} = \sum_{i=1}^{N-k} p_{i,i+k}$ and $\bar{s}_{kN} = \sum_{i=k+1}^{N} p_{i,i-k}$. First, there is a simple recursion formula for the sums $s_{kN}$ that follows from the recursion formula for the $p_{ij}$'s (i.e., Theorem 2). We claim that
$$s_{k+1,N+1} = s_{k-1,N} + 2s_{kN} + s_{k+1,N}.$$
The proof amounts to computing the sum of hook-sum matrices $H_{ij}^N$ whose vertices all lie on the $k$-th diagonal above the main diagonal (i.e., $j - i = k$). Each entry is a sum of +1's and -1's, but the key point is that there are always more +1's than -1's, so in the sum we are left with only 1's along the $k$-th diagonal, 2's along the $(k+1)$-th, and 1's along the $(k+2)$-th. Finally, some simple renumbering of the diagonals has to be done when we move from an $N \times N$ matrix to $(N+1) \times (N+1)$.

Although the recursion formula for the diagonal sums looks simple, it actually is quite important. This is the only step of the solution of Sun Bin's problem that requires a recursion on $N$, and the only step that uses the hook-sum map. The precision of this recursion is what makes it possible to get such a precise calculation of $k^*(N)$.

There is also a recursion formula for $\bar{s}_{k,N+1}$, but it is not as simple. The reason is that in the corresponding sum of hook-sum matrices, when the vertices all lie on the $k$-th diagonal *below* the main diagonal, some of the -1's do not cancel. However, we can work around this because by Lemma 1.3,

$$\bar{s}_{kN} = (N-k)\binom{2N}{N} - s_{kN}.$$

So we can simply compute $\bar{s}_{kN}$ from $s_{kN}$.

Now remember that we are looking for the *maximum* diagonal sum below the main diagonal, because (assuming the No Gaps Conjecture) the optimal permutation matrix will have 1's along one of these diagonals. Eventually we will also have to take account of the 1's above the main diagonal as well, but for now we will ignore those "nuisance" terms.

The maximum sum will occur at a value of $k$ such that $\Delta\bar{s}_{kN} \equiv \bar{s}_{kN} - \bar{s}_{k-1,N} \geq 0$ but $\Delta\bar{s}_{k+1,N} < 0$. From numerical examples it is quickly apparent that the differences $\Delta\bar{s}_{1N}, \Delta\bar{s}_{2N}, \ldots$ form a strictly decreasing sequence, and hence there is a unique value of $k$, $k_{\max}$, with this property. However, we have not yet proved that the sequence $\{\Delta\bar{s}_{kN}\}$ is decreasing. This will follow from a remarkably simple closed-form expression for the terms of the sequence $\{\Delta s_{kN}\}$:

$$\Delta s_{kN} = -\binom{2N}{0} - \binom{2N}{1} - \ldots - \binom{2N}{N-k}.$$

The explicit formula for $\Delta s_{kN}$ can be proved by induction on $N$. First of all, for $N=2$, for which $P^2 = \begin{bmatrix} 3 & 1 \\ 5 & 3 \end{bmatrix}$, we have $s_{0N} = 6$ and $s_{1N} = 1$, so $\Delta s_{1N} = 1 - 6 = -5 = -\binom{4}{0} - \binom{4}{1}$. Similarly, $\Delta s_{2N} = -1 = -\binom{4}{0}$. The inductive step from $N$ to $N+1$ is very easy and is left to the reader; the key point is that the entries in the $(2N+2)$-th row of Pascal's triangle satisfy exactly the same recursion that the $s_{kN}$'s do, namely

$$\binom{2N+2}{k+1} = \binom{2N}{k-1} + 2\binom{2N}{k} + \binom{2N}{k+1}.$$

From the formula for $s_{kN}$ and the symmetry relation between $s_{kN}$ and $\bar{s}_{kN}$, it follows immediately that

$$\Delta\bar{s}_{kN} = \binom{2N}{0} + \binom{2N}{1} + \ldots + \binom{2N}{N-k} - \binom{2N}{N}.$$

We pause here to point out that because $s_{0N}$ is known and the first differences $\Delta s_{kN}$ and $\Delta\bar{s}_{kN}$ are known, it is an easy matter to compute $s_{kN}$ and $\bar{s}_{kN}$. This completes the proof of Lemma 1.7 (postponed from Section 1).

However, it is the first differences that are crucial for solving the maximization problem, and to which we will turn our attention for the rest of the proof. As mentioned above, it is useful to have a "center-out" formula for the first differences $\Delta\bar{s}_{kN}$. It is well-known that $\sum_{j=0}^{2N}\binom{2N}{j} = 4^N$. Hence,

$$\binom{2N}{0} + \binom{2N}{1} + \ldots + \binom{2N}{N-1} + \frac{1}{2}\binom{2N}{N} = 2 \cdot 4^{N-1} \text{ and}$$

$$\Delta\bar{s}_{kN} = 2 \cdot 4^{N-1} - \frac{3}{2}\binom{2N}{N} - \binom{2N}{N-1} - \ldots - \binom{2N}{N-k+1}.$$

This makes the decreasing property of $\Delta \bar{s}_{kN}$ exceedingly obvious, since one additional negative term is added each time $k$ is incremented by 1. It is also easy to show that for $N \geq 3$, $\Delta \bar{s}_{1N} > 0$. Thus for all $N \geq 3$, there is in fact a unique integer,
$$k(N) = \sup\{k : \Delta \bar{s}_{kN} \geq 0\} = \inf\{k : \Delta \bar{s}_{kN} < 0\} - 1,$$
which satisfies the conditions necessary for a maximum diagonal sum.

Finally, we address the issue of the "nuisance terms." By assuming the No Gaps Conjecture, we are assuming that the optimal permutation is one of the permutations
$$\pi_k = \{N, N-1, ..., N-k+1, 1, 2, ..., N-k\}.$$
The objective function $F(\pi_k)$ is the sum of two terms: the sum of the $p_{ij}$'s along the $k$-th diagonal below the main diagonal (i.e., $\bar{s}_{kN}$), and the sum of the first $k$ terms along the anti-main diagonal. Thus, using Lemma 1.2,
$$F(\pi_k) - F(\pi_{k-1}) = \Delta \bar{s}_{kN} + p_{k,N+1-k} = \Delta \bar{s}_{kN} + \binom{N}{k-1}^2.$$
The last term here accounts for the extra "nuisance term" in the definition of $k^*(N)$. In addition it is necessary to prove that the sequence $\{F(\pi_k) - F(\pi_{k-1})\}$ is strictly decreasing and hence that there is a single value of $k$, $k^*(N)$, where it crosses from positive to negative. As it turns out, this is not quite straightforward, because the first term ($\Delta \bar{s}_{kN}$) is decreasing but the nuisance term $\binom{N}{k-1}^2$ is increasing. What we can show instead is that in the relevant range where the crossover can occur, say for $k \leq N/2$, the nuisance term is dominated by the last term in the sum for $\Delta \bar{s}_{kN}$, and hence the crossover occurs at most one value of $k$ later and is unique within this range. The details are left to the reader. □

While Theorem 1 in principle solves the problem of finding the optimal permutation (assuming the No Gaps Conjecture for $N \leq 10^7$ and unconditionally for $N > 10^7$), it would be nice to have an estimate of the rate of growth of $k(N)$ and $k^*(N)$. We will show that $k(N)$ is bounded above and below by constants times $\sqrt{N \ln N}$, for sufficiently large $N$. Because $k^*(N)$ differs from $k(N)$ by at most one, the same estimate holds for it.

To estimate $k(N)$, we want to answer the question: How many binomial coefficients $\binom{2N}{j}$ do we need to add up to make
$$\binom{2N}{0} + \binom{2N}{1} + ... + \binom{2N}{N-k} \geq \binom{2N}{N}?$$
The left-hand side of this equation is simply the cumulative distribution function (CDF) for the binomial distribution, which is known to approach the CDF for the normal distribution, by the Central Limit Theorem. Hence what we need is an estimate of the rate of convergence in the CLT. This question has been well-studied in the literature, and the state-of-the-art estimates go by the name of the Berry-Esséen Theorem.

There is no doubt that an asymptotic estimate for $k(N)$ can be derived from the Berry-Esséen Theorem. However, this theorem is not ideally suited for the purpose, because it is a theorem about a distribution function but the question above involves a "point mass," the central binomial coefficient. To put it another way, we don't want an absolute measure of the error in the normal approximation, but a relative measure, scaled to $\binom{2N}{N}$.

Very much to my surprise, I was able to find an elementary relative estimate (using nothing more sophisticated than Taylor's Theorem) that *improves* upon Berry-Esséen *in the region of greatest interest to us*. Specifically, Berry-Esséen gives an estimate of the form -ε < error < ε, but I was able to prove a one-way estimate, 0 < error < ε, for a surprisingly large interval of $k$. Of course this estimate cannot hold for all $k$ because that would violate the definition of a probability distribution, but the values of $k$ for which it fails are very far from $k(N)$ and hence do not enter into our estimate of $k(N)$.

I don't know if this "one-sided" behavior has been noted previously in the literature, so the following lemma may be of independent interest.

**Lemma 4.5. (One-Sided Convergence in the Central Limit Theorem).** If $N \geq 8$ and $3 \leq k \leq \dfrac{N}{2} - 1$, then

$$\binom{2N}{N-k} < e^{-k^2/N} \binom{2N}{N}.$$

**Proof.** First we note that

$$\frac{\binom{2N}{N-k}}{\binom{2N}{N}} = \frac{N(N-1)\cdots(N-k+1)}{(N+1)(N+2)\cdots(N+k)}.$$

We apply the arithmetic-geometric mean inequality and then bound the right-hand side by an integral:

$$\left[\frac{N(N-1)\cdots(N-k+1)}{(N+1)(N+2)\cdots(N+k)}\right]^{1/k} \leq \frac{1}{k}\left[\frac{N}{N+k} + \frac{N-1}{N+k-1} + \ldots + \frac{N-k+1}{N+1}\right]$$

$$\leq \int_{(N-k+1)/k}^{(N+1)/k} \frac{x}{x+1}dx = 1 - \ln\left(1 + \frac{k}{N-k+1}\right) < 1 - \frac{k}{N-k+1} + \frac{1}{2}\left(\frac{k}{N-k+1}\right)^2,$$

where the last inequality follows from Taylor's Theorem.

Next we note that

$$-\frac{k}{N-k+1} + \frac{1}{2}\left(\frac{k}{N=k+1}\right)^2 + \frac{k}{N} = \frac{k}{N(N-k+1)^2}[k^2 - (2+N/2)k + N + 1].$$

The quadratic polynomial on the right-hand side is less than or equal to zero for any value of $k$ between its two roots, $1 + \dfrac{N}{4} \pm \dfrac{N}{4}\sqrt{1 - 8/N}$. (Here we use the assumption that $N \geq 8$. For smaller values of $N$ the roots are not real.) Note that the two roots are equally spaced about $1+N/4$. Also notice that the polynomial in braces is always less than or equal to zero for $k = 3$ (again assuming $N \geq 8$). Thus the smaller root of the quadratic must be less than or equal to 3, and this means by symmetry that the larger root must be greater than or equal to $N/2 - 1$. Putting it all together, we have

$$-\frac{k}{N-k+1} + \frac{1}{2}\left(\frac{k}{N-k+1}\right)^2 \leq -\frac{k}{N}$$

for $N \geq 8$ and $3 \leq k \leq \frac{N}{2} - 1$. Hence, for the same range of $N$ and $k$,

$$\frac{N(N-1)\cdots(N-k+1)}{(N+1)(N+2)\cdots(N+k)} \leq (1-k/N)^k = e^{k\ln(1-k/N)} < e^{-k^2/N}. \qquad \square$$

To get a lower bound on $k(N)$, we want to ask for which values of k we can be sure that $\Delta \bar{s}_{kN} \geq 0$. We'll start with a rough estimate and then refine it to a more accurate estimate.

**Proposition 4.6. (Lower Bound)** For all $N, k(N) \geq \lfloor (\sqrt{\pi N} - 1)/2 \rfloor$. For all $N \geq 400$, $k(N) \geq \sqrt{N \ln N / 4}$. Asymptotically, for any $\varepsilon > 0$ there exists $N_\varepsilon > 0$ such that for all $N$ greater than $N_\varepsilon$, $k(N) \geq (1-\varepsilon)\sqrt{N \ln N / 2}$.

**Proof.** We start from the alternative formula for $\Delta \bar{s}_{kN}$:

$$\Delta \bar{s}_{kN} = 2 \cdot 4^{N-1} - \left[ \frac{3}{2} + \frac{N}{N+1} + \frac{N(N-1)}{(N+1)(N+2)} + \ldots + \frac{N\cdots(N-k+2)}{(N+1)\cdots(N+k)} \right] \binom{2N}{N}.$$

The simplest, most naïve estimate, not even using Lemma 4.5, is to notice that each of the quotients in braces, after $3/2$, is less than 1. Hence

$$\Delta \bar{s}_{kN} > 2 \cdot 4^{N-1} - (k + \frac{1}{2})\binom{2N}{N}.$$

A standard estimate on the central binomial coefficient says that $\binom{2N}{N} < 4^N / \sqrt{\pi N}$. Hence

$$\Delta \bar{s}_{kN} > 4^N \left[ \frac{1}{2} - (k + \frac{1}{2})/\sqrt{\pi N} \right] > 0$$

provided $k < (\sqrt{\pi N} - 1)/2$. Because $k(N) = \sup\{k \in Z^+ : \Delta \bar{s}_{kN} \geq 0\}$, the first conclusion of the theorem follows.

This naïve estimate is not bad for small values of $N$. For example, for $N = 8$ it tells us that $k(8) \geq 2$ (so that in a game of One Round War with 8 cards, we should throw at least 2 tricks). In fact, $k(8) = 2$, so the inequality is sharp. Likewise, $k(16) \geq 3$ according to the naïve estimate, and indeed it turns out that $k(16) = 3$, so the inequality is sharp. However, there are other cases where it is not sharp. For example, the naïve estimate gives $k(13) \geq 2$, when in fact $k(13) = 3$. As $N$ gets large, the naïve estimate becomes less sharp, and we would like to get a better estimate.

To get a better estimate, we approximate the sum of binomial coefficients by an integral, specifically $\int_0^k e^{-x^2/N} dx$. The presence of the "3/2" suggests that we should use the Trapezoid Rule expansion, which will borrow ½ from that first term. We use Lemma 4.5 to make sure that our approximation is a lower bound. Of course, Lemma 4.5 does not apply to the first two quotients after the "3/2" term, because of the requirement that k ≥ 3 in that lemma. Thus, the first two terms have to be estimated by hand. Next we write the integral as $\int_0^\infty e^{-x^2/N} dx - \int_k^\infty e^{-x^2/N} dx$. We note that the first integral can be evaluated exactly and, in fact, it almost exactly cancels the $2 \cdot 4^N$ term on account of the asymptotic estimate $\binom{2N}{N}/4^N \sim 1/\sqrt{\pi N}$.

When we put all these ideas together we obtain that

$$\Delta \bar{s}_{k+1,N} \geq \binom{2N}{N}\left[-1 + \int_k^\infty e^{-x^2/N}dx - \varepsilon(k,N)\right],$$

where $\varepsilon(k, N)$ represents all the error terms—from the Trapezoid Rule, from the first two terms of the sum, and from the asymptotic estimate for the central binomial coefficient. Of course, for a precise estimate the error terms have to be kept track of carefully, but they are all small compared to the first two terms. The next step is to integrate by parts twice on the right-hand side and notice that the remaining integral is positive, which leads to the estimate

$$\int_k^\infty e^{-x^2/N}dx > e^{-k^2/N}\left(\frac{N}{2k} - \frac{N^2}{4k^3}\right).$$

Plugging this into the inequality above and doing a careful analysis of the error terms leads to the following estimate:

$$\Delta \bar{s}_{k+1,N} \geq \binom{2N}{N}\left[-1 + e^{-k^2/N}\left(\frac{N}{2k} - \frac{N^2}{4k^3} - \frac{1}{2}\right) - \frac{k}{3e^{3/2}N}\right].$$

Now let $k = \left\lfloor \frac{1}{2}\sqrt{N \ln N} \right\rfloor$. If $N \geq 400$, I claim that $\Delta \bar{s}_{k+1,N} > 0$. Because $\Delta \bar{s}_{k+1,N}$ is decreasing as a function of $k$, it will follow that $\Delta \bar{s}_{kN} > 0$ for all $k \leq \frac{1}{2}\sqrt{N \ln N}$, so that in fact $k(N) > \frac{1}{2}\sqrt{N \ln N}$.

The key point to realize is that for $k = \left\lfloor \frac{1}{2}\sqrt{N \ln N} \right\rfloor$, $k^2/N \leq \ln N/4$, so $e^{-k^2/N} \geq N^{-1/4}$.

Meanwhile, the positive term in parentheses is $N/2k$, which is on the order of $N^{1/2}/\ln N$. So the product is on the order of $N^{1/4}/\ln N$, as $N$ approaches infinity. The largest of the negative terms, $e^{-k^2/N}\left(\frac{N^2}{4k^3}\right)$, is on the order of $N^{1/4}/(\ln N)^2$, and the other negative terms are all either constant or decreasing. Hence it is obvious that for large enough $N$, the positive term will be dominant and $\Delta \bar{s}_{k+1,N} > 0$. Verifying that 400 is "large enough" is just a matter of plugging the numbers into a calculator.

As for the asymptotic estimate, the argument is just the same. If we take $k = c\sqrt{N \ln N}$, where $c^2 = 1/2 - \delta$ for some small positive $\delta$, then $e^{-k^2/N}\left(\frac{N}{2k}\right)$ is of order $N^\delta/\ln N$, while all the negative terms are of lower order. Hence for sufficiently large $N$, the positive term must dominate and $\Delta \bar{s}_{k+1,N} > 0$, from which it follows that $k(N) > c\sqrt{N \ln N}$. $\square$

**Remark 1.** Some readers might wonder whether the Trapezoid Rule was really necessary in this proof. The answer is no, but it made the error terms smaller and hence improved the estimate of when the $N^{1/4}/\ln N$ term begins to dominate the others.

**Remark 2.** As mentioned earlier, Proposition 4.6 in fact holds for $N \geq 91$, as shown by computer calculations. The empirical evidence also suggests that the rate at which $k(N)/\sqrt{N \ln N}$ approaches its eventual limit of $1/\sqrt{2}$ is extraordinarily slow. Even at $N = 500$, this quotient is hovering in a range between 0.54 and 0.56, and there is little to suggest that it will eventually converge to 0.707… However, the data do reveal one interesting phenomenon. Aside from the jumps caused by the fact that $k(N)$ is an integer-valued function, the overall trend of $k(N)/\sqrt{N \ln N}$ seems to be monotone increasing, and this suggests that it approaches its limit strictly from below. This observation is confirmed by the following proposition.

**Proposition 4.7 (Upper Bound).** For all $N \geq 30$, $k(N) < \sqrt{N \ln N / 2}$.

**Proof.** In this case we will make the argument easier to follow by using previously published estimates. According to Theorem 5.3.2 in [LPV],

$$\sum_{i=0}^{N-k-1} \binom{2N}{i} \leq 2^{2N-1} \exp\left(\frac{-k^2}{N+k}\right).$$

If we set $k = \lceil \sqrt{N \ln N / 2} \rceil$, the right-hand side is less than $\frac{1}{2} \cdot 4^N \cdot N^{-N/2(N+k)}$. To prove the proposition, we simply need to determine when this sum is less than $\binom{2N}{N}$. Here it is useful to apply the lower bound $\binom{2N}{N} > \frac{4^N}{\sqrt{\pi N}}\left(1 - \frac{1}{8N}\right)$, which can be found in the Wikipedia entry on "central binomial coefficient." From some simple calculations it follows that we need

$$\frac{\sqrt{\pi}}{2} \frac{N^{\sqrt{\ln N / 8N}}}{1 - (1/8N)} < 1.$$ The numerator has its maximum value at $N = e^3$ and is decreasing thereafter. So it suffices to find a value of $N$ for which the numerator is slightly less than $2/\sqrt{\pi}$. (The denominator will be very close to 1 by this point.) Trial and error on a calculator show that $N = 10{,}000$ works, hence the proposition is true for all $N \geq 10{,}000$.

To improve the estimate to $N \geq 30$, we can use the strategy we used for Proposition 4.6 (i.e., employ the normal approximation to the CDF of the binomial distribution and do careful error estimates). I will not show the details, which are similar to the proof of Proposition 4.6 and extremely tedious, but will leave it to motivated readers to work out for themselves. Incidentally, computer calculations show that Proposition 7 is in fact true for $3 \leq N \leq 30$ as well. □

**Remark 1.** It seems strange that the book [LPV] gives only an upper bound for the sum of the binomial coefficients, which was just what we needed for Proposition 4.7. It would be nice if there were also a published lower bound, but I couldn't find one in the literature that was strong enough to give an explicit cutoff for $N$. So for Proposition 4.6 I was forced to improvise.

**Remark 2.** As pointed out in section 1, the state-of-the-art general-purpose algorithm for solving linear assignment problems runs in time proportional to $O(N^3)$. By contrast, Theorem 4.1 can be viewed as a special-purpose algorithm for solving Sun Bin's problem only. If we use the "center-out" sum, the algorithm consists of adding at most $(N \ln N/2)^{1/2}$ numbers together and comparing them with $\binom{2N}{N}$.

Thus, if we were given an oracle to produce Pascal's triangle on demand, then our special-purpose algorithm would run in time $O((N \ln N)^{1/2})$, or sublinear time. Of course, it is unrealistic to assume the existence of such an oracle. But even without it, the $(2N)$-th row of Pascal's triangle can be generated in at most $O(N^2)$ steps, and thus the special-purpose Sun Bin Algorithm is faster than the general-purpose Hungarian Algorithm.

**Bibliography.**

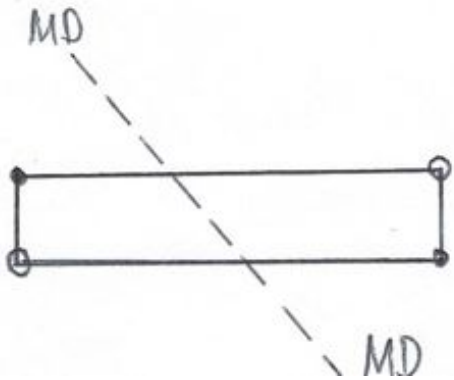

**Figure 1.1.** A "bad" rectangle. It is impossible (?) to find a general rule that states which permutation has a higher value of the objective function – the hollow dots or the filled dots.

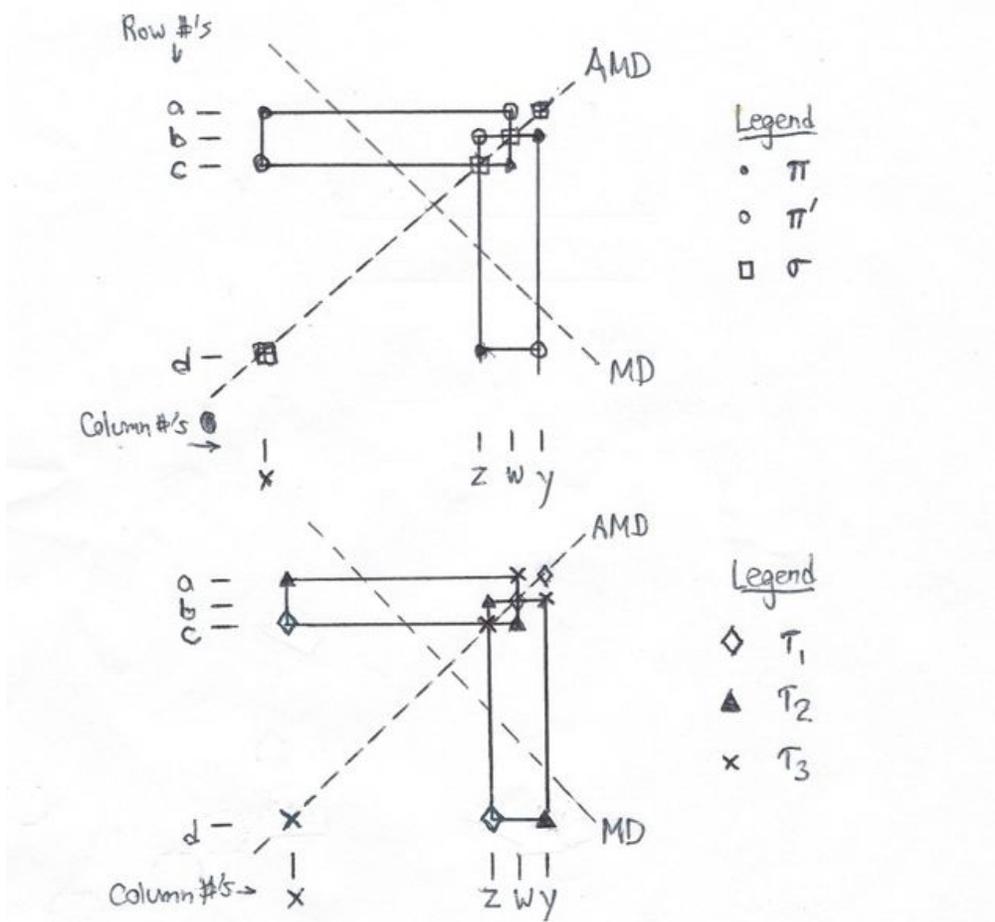

**Figure 2.1.** The six permutations π, π′, σ, τ₁, τ₂, τ₃ in the simplest case, where π and π′ differ by two 2-cycles. AMD = anti-main diagonal, MD = main diagonal. Note that π′ is the reflection of π through the main diagonal, while σ, τ₁, τ₂, τ₃ are all symmetric.

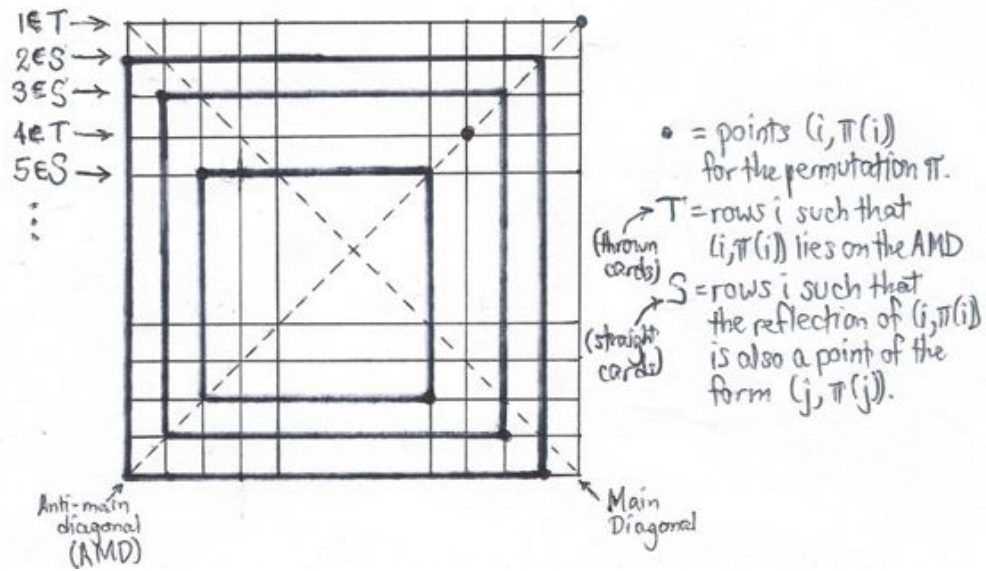

**Figure 2.2.** The "outside-in" inductive argument for the proof of the Symmetry Lemma.

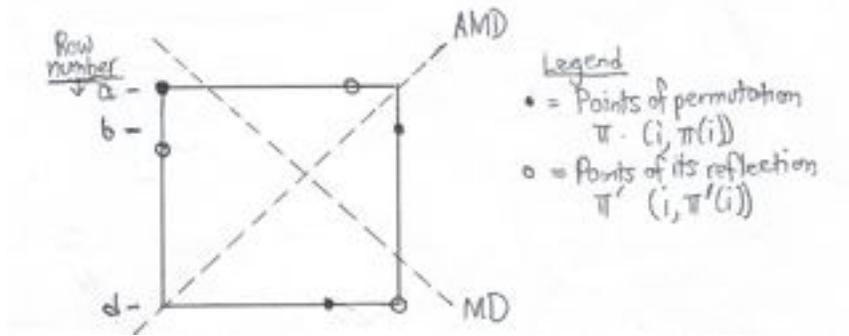

**Figure 2.3.** Known information about π and π′ after application of Inductive Hypothesis. Row *a* is the first row where π is believed to be non-symmetric. Note that there is already some resemblance to Figure 2.1.

|  | $\pi$ | $\pi'$ | $\sigma$ | $\tau_1$ | $\tau_2$ | $\tau_3$ |
|---|---|---|---|---|---|---|
| First r-cycle | $\pi(b)=pa$ $\pi(a_1)=x_1$ $\vdots$ $\pi(a_{r-2})=x_{r-2}$ $\pi(d)=x_{r-1}$ | $\pi'(b)=x_1$ $\pi'(a_1)=x_2$ $\vdots$ $\pi'(a_{r-2})=x_{r-1}$ $\pi'(d)=pa$ | $\sigma(b)=pb$ $\sigma(a_1)=pa_1$ $\vdots$ $\sigma(a_{r-2})=pa_{r-2}$ $\sigma(d)=pd$ | $\tau_1(b)=pb$ $\tau_1(a_1)=x_1$ $\vdots$ $\tau_1(a_{r-2})=x_{r-2}$ $\tau_1(d)=x_{r-1}$ | $\tau_2(b)=x_1$ $\tau_2(a_1)=x_2$ $\vdots$ $\tau_2(a_{r-2})=x_{r-1}$ $\tau_2(d)=pa$ | $\tau_3(b)=pa$ $\tau_3(a_1)=pa_1$ $\vdots$ $\tau_3(a_{r-2})=pa_{r-2}$ $\tau_3(d)=pd$ |
| Second r-cycle | $\pi(a)=pd$ $\pi(px_1)=pb$ $\pi(px_2)=pa_1$ $\vdots$ $\pi(px_{r-1})=pa_{r-2}$ | $\pi'(a)=pb$ $\pi'(px_1)=pa_1$ $\pi'(px_2)=pa_2$ $\vdots$ $\pi'(px_{r-1})=pd$ | $\sigma(a)=pa$ $\sigma(px_1)=x_1$ $\sigma(px_2)=x_2$ $\vdots$ $\sigma(px_{r-1})=x_{r-1}$ | $\tau_1(a)=pa$ $\tau_1(px_1)=pa_1$ $\tau_1(px_2)=pa_2$ $\vdots$ $\tau_1(px_{r-1})=pd$ | $\tau_2(a)=pd$ $\tau_2(px_1)=pb$ $\tau_2(px_2)=pa_1$ $\vdots$ $\tau_2(px_{r-1})=pa_{r-2}$ | $\tau_3(a)=pb$ $\tau_3(px_1)=x_1$ $\tau_3(px_2)=x_2$ $\vdots$ $\tau_3(px_{r-1})=x_{r-1}$ |

**Figure 2.4.** Case 1 of Symmetry Lemma: "Splicing together" $\pi$ and $\pi'$. Colors are intended to make it easy to see the splicing. The permutation $\tau_2$ is obtained by splicing the first half of $\pi'$ with the second half of $\pi$; $\tau_1$ by splicing the first half of $\pi$ with the second half of $\pi'$ and "transfecting" the result with row $a$ and row $b$ of $\sigma$. The permutation $\tau_3$ differs from $\sigma$ only in rows $a$ and $b$. In case $r = 2$, all equations involving $a_i$ are vacuous, and only $x_1$ is defined (of the $x$'s). In this case the permutations are as shown in Figure 2.1.

The reader should check that the domain and range of each map is the same, and each map is bijective. In addition, the reader should check that in each row, the images of $\tau_1$, $\tau_2$, and $\tau_3$ are a permutation of the images of $\pi$, $\pi'$, and $\sigma$.

**Figure 2.5.** "Splicing" together $\pi$ and $\pi'$ in Case 2. Here $r$ must be odd, and $m = (r+1)/2$. Colors are again intended to make it easy to see the splicing. The splicing procedure is slightly different from Case 1; one of the main differences is that both $\tau_1$ and $\tau_2$ get contributions from $\sigma$ (red).

The same comments ("the reader should check that…") hold for this figure as for Figure 2.4.

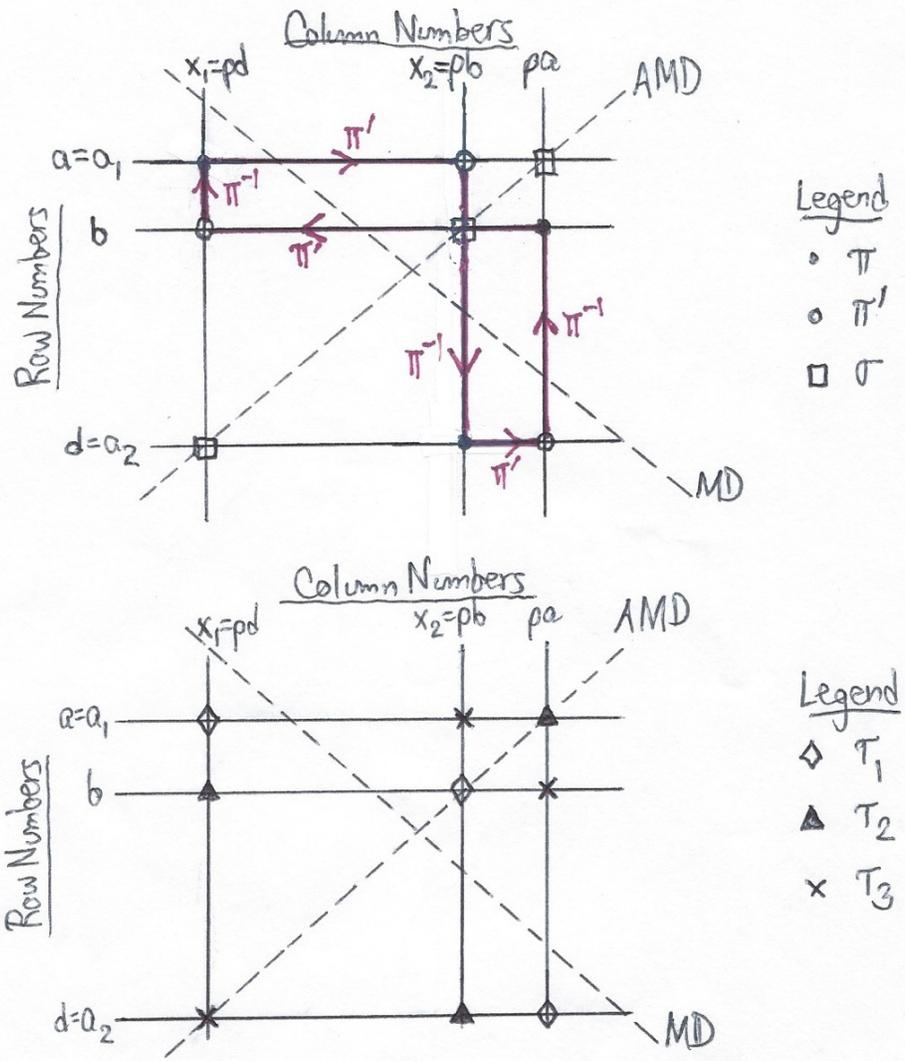

**Figure 2.6.** Case 2, illustrated for the sub-case $r = 3$ (where certain equations in Figure 2.5 become vacuous). The 3-cycle $\pi^{-1}\pi'$ is highlighted in red. Note that in this subcase, all six bijections from $\{b, a_1, a_2\}$ to $\{\rho a, x_1, x_2\}$ are used in the proof.

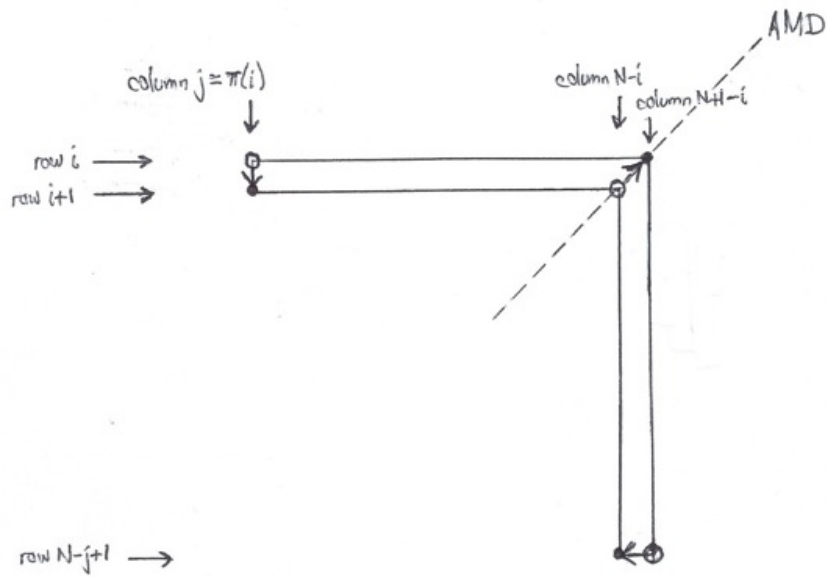

**Figure 3.1.** Here $\pi$ is assumed to be a permutation with a "rogue element," i.e., it has a "1" in the $(i+1)$-th position on the AMD and a "0" in the $i$-th position. Hollow dots represent the permutation $\pi$. We can get a new permutation (filled dots) by sliding the rogue element one step to the northeast, while sliding the "1" in row $i$ down one row and the "1" in column $(N+1-i)$ one column to the left. As usual, AMD denotes the anti-main diagonal.

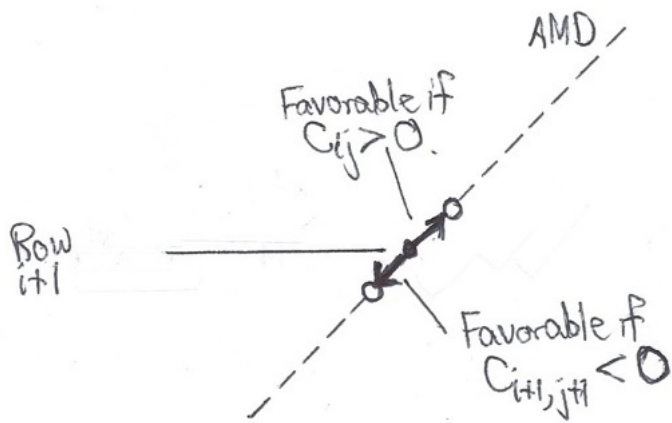

**Figure 3.2.** Here $C_{ij}$ denotes the change to the objective function resulting from the "bead sliding" move in Figure 3.1. The simplest case is that of an isolated rogue element in row $(i + 1)$. If $C_{ij} > 0$, we can increase the objective function by sliding it one step to the northeast. If $C_{i+1,j+1} < 0$, we can increase the objective function by sliding it one step to the southwest. However, if $C_{i+1,j+1} > 0 > C_{ij}$, then the bead is stuck. The fundamental challenge in proving the No Gaps Theorem is either to deal with stuck beads or show that they cannot happen.

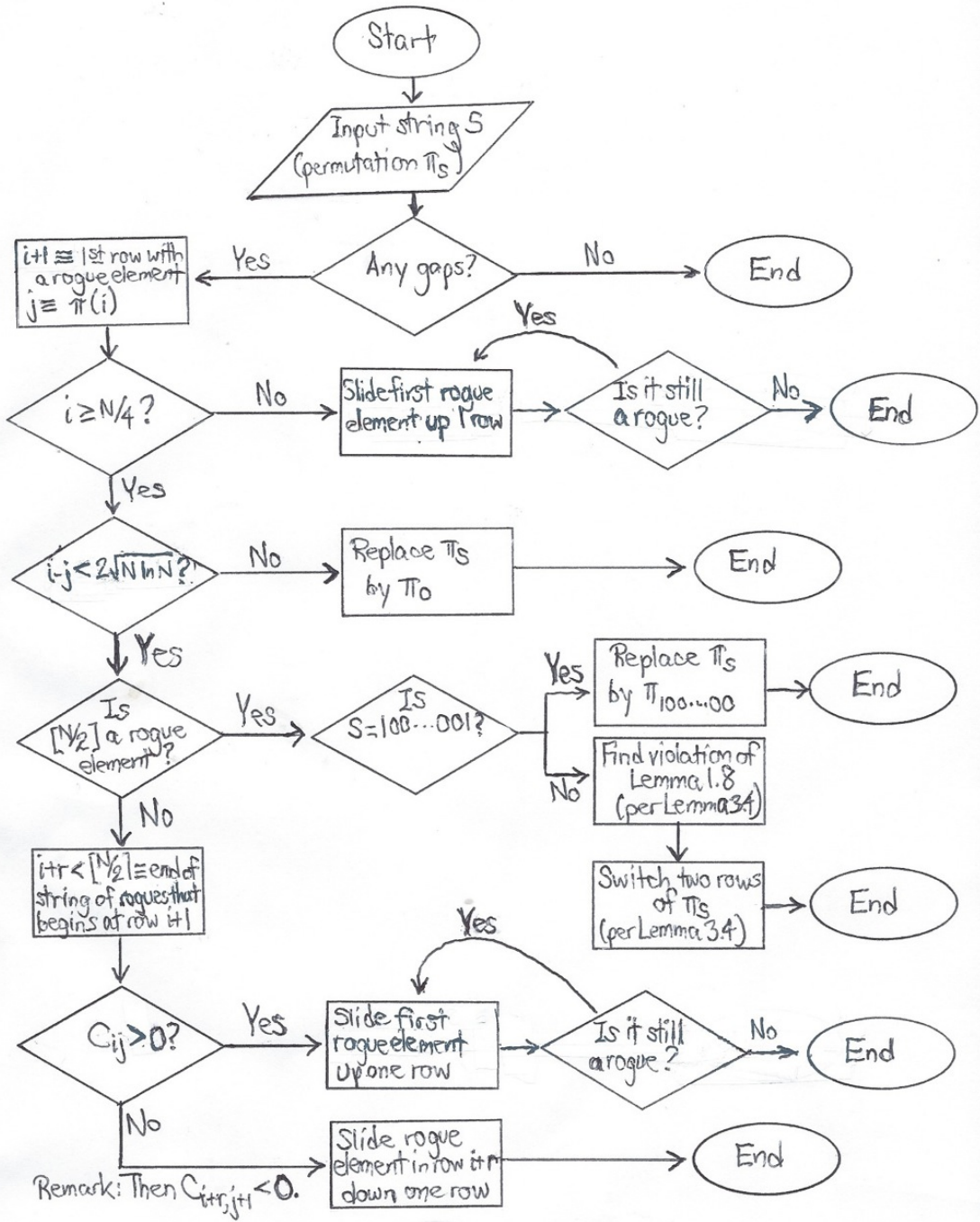

**Figure 3.3.** Flow chart for the proof of the No Gaps Theorem.

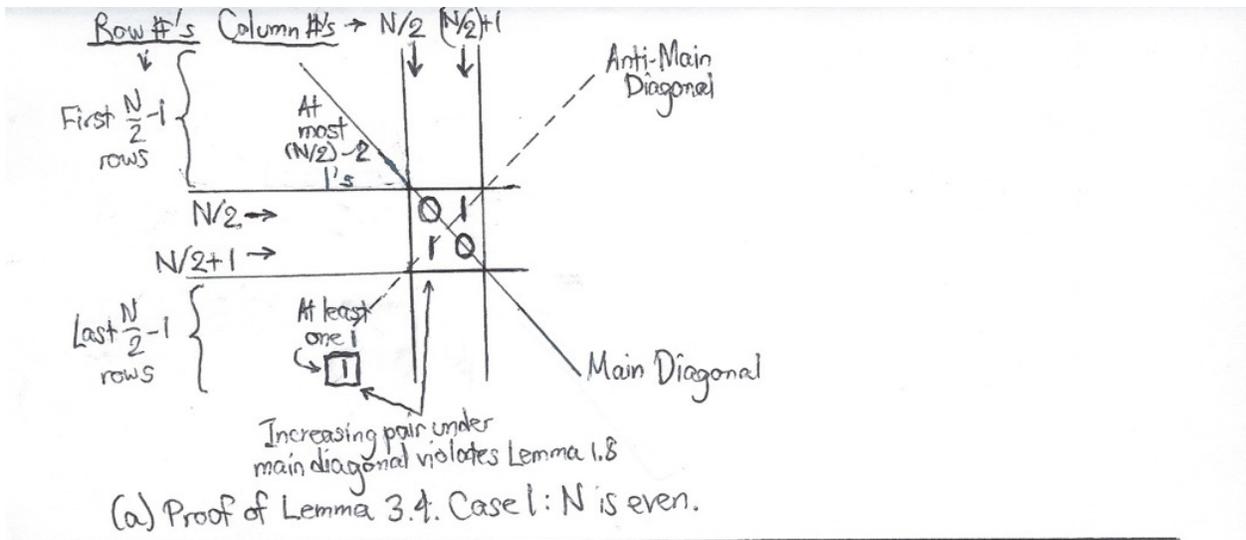

(a) Proof of Lemma 3.4. Case 1: N is even.

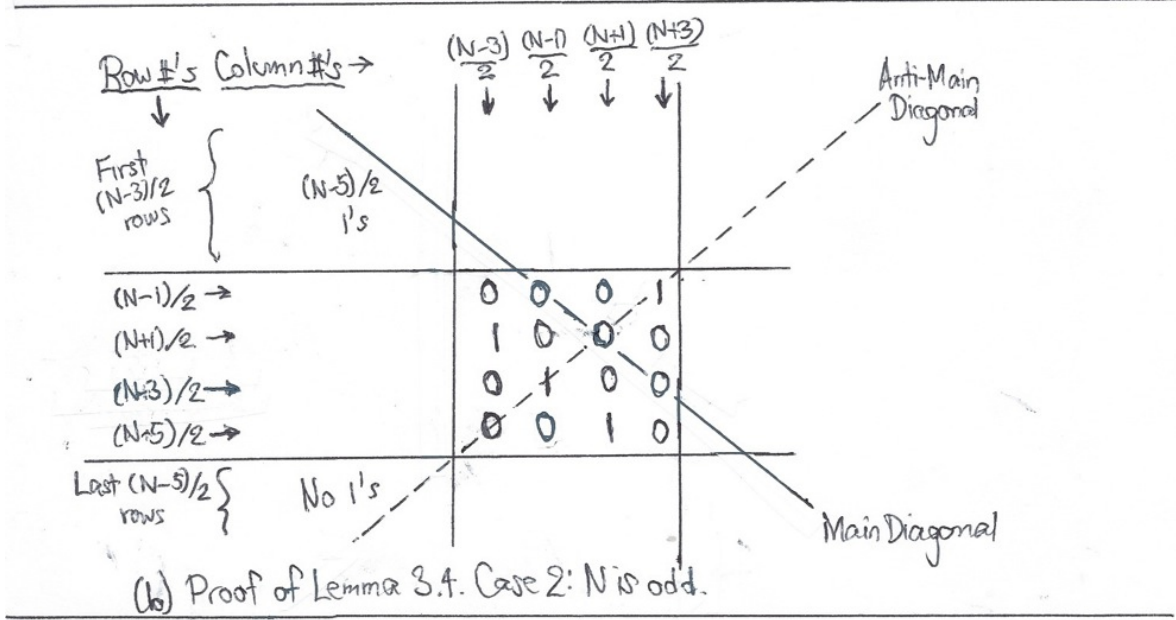

(b) Proof of Lemma 3.4. Case 2: N is odd.

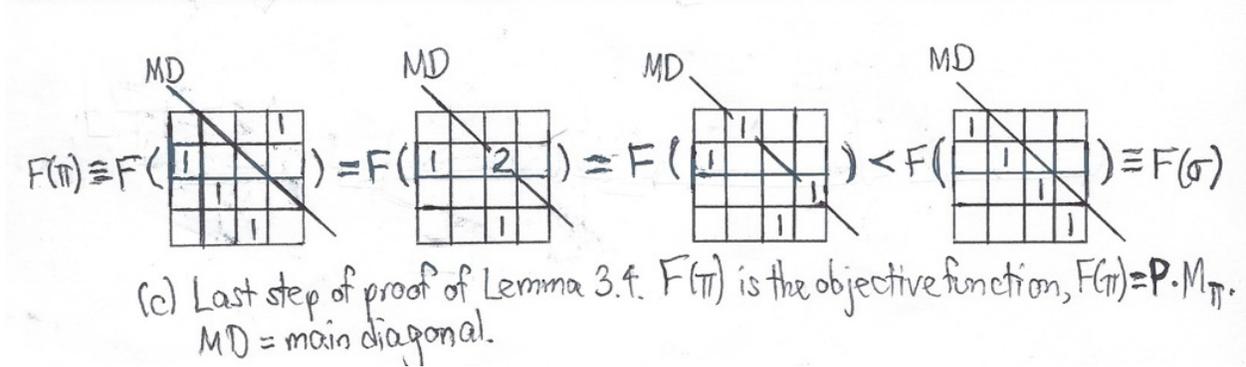

(c) Last step of proof of Lemma 3.4. $F(\pi)$ is the objective function, $F(\pi) = P \cdot M_\pi$. MD = main diagonal.

Figure 3.4. Visual proof of Lemma 3.4.